\documentclass[12pt]{amsart}
\usepackage{amscd,amssymb,verbatim}

\usepackage[T1]{fontenc}
\usepackage[utf8]{inputenc}
\usepackage{lmodern}
\usepackage{amsmath,amssymb,amsthm}
\usepackage[margin=29mm]{geometry}
\usepackage{enumitem}
\usepackage{microtype}
\usepackage[hidelinks]{hyperref}
\usepackage{xurl}

\newtheorem{theorem}{Theorem}[section]
\newtheorem*{introtheorem}{Theorem~\ref{thm:abstract}}
\newtheorem*{introoperator}{Theorem~\ref{thm:main}}
\newtheorem{proposition}[theorem]{Proposition}
\newtheorem{corollary}[theorem]{Corollary}
\newtheorem{lemma}[theorem]{Lemma}
\newtheorem{definition}[theorem]{Definition}

\newtheorem{question}{Question}
\theoremstyle{remark}
\newtheorem{remark}[theorem]{Remark}

\newcommand{\R}{\mathbb R}
\newcommand{\N}{\mathbb N}
\newcommand{\Q}{\mathbb Q}
\newcommand{\clspan}{\overline{\operatorname{span}}}
\newcommand{\norm}[1]{\left\lVert#1\right\rVert}

\DeclareMathOperator{\supp}{supp}
\DeclareMathOperator{\dens}{dens}

\title{A description of the space with no separable quotient presented by OpenAI}
\author{Spiros A. Argyros}
\address{National Technical University of Athens, Faculty of Applied Sciences, Department of Mathematics, Zografou Campus, 157 80, Athens, Greece.}
\email{sargyros@math.ntua.gr}

\begin{document}
\begin{abstract}
OpenAI recently presented a remarkable solution to the longstanding
separable quotient problem in Banach space theory. The result establishes
its independence from ZFC, relative to the consistency of ZFC with a
measurable cardinal. Under the continuum hypothesis, the paper constructs
a nonseparable Banach space $X$ with no separable infinite-dimensional
quotient. The present paper offers an alternative exposition of the
mechanism leading to $X$ and its properties. Our aim is to make the
methods and results more accessible and to offer a clearer understanding
of this exotic space.
\end{abstract}

\maketitle


\section*{Introduction}

The separable quotient problem asks whether every infinite-dimensional
Banach space admits a separable infinite-dimensional quotient. Equivalently,
given an infinite-dimensional Banach space $X$, does there exist a closed
subspace $Y$ of $X$ such that $X/Y$ is separable and infinite dimensional?
One of the attractive features of this problem is that its formulation is
accessible to anyone who has completed a first course in functional
analysis. Every infinite-dimensional Banach space contains a separable
infinite-dimensional closed subspace, obtained by taking the closed linear
span of a countable linearly independent set. The corresponding question
for quotients has proved much more difficult.

Traditionally associated with Banach and Mazur, the problem has occupied
a central place in the structural theory of Banach spaces. Among the
positive results preceding the work discussed here, Johnson and Rosenthal
\cite{JR72} proved in 1972 that every separable infinite-dimensional
Banach space has an infinite-dimensional quotient with a Schauder basis.
In 2008, Argyros, Dodos and Kanellopoulos \cite{ADK08} established that
every infinite-dimensional dual Banach space has a separable
infinite-dimensional quotient.

The nontrivial case of the separable quotient problem concerns
nonseparable spaces, since a separable infinite-dimensional space is
already a quotient of itself. In the separable setting, a related question
asks whether there exists a separable infinite-dimensional Banach space
such that every infinite-dimensional quotient has a nonseparable dual.
This question, posed by Johnson and Rosenthal \cite{JR72}, has a positive
answer. Such a space is supplied by the construction of Argyros,
Arvanitakis and Tolias in \cite[Section~6]{AAT06}; the quotient property
is explicitly established by Argyros and Beanland in
\cite[Section~5, Corollary~22]{AB13}. Thus nonseparability of the dual may
persist in every infinite-dimensional quotient of a separable space.

In the preprint \emph{Relative independence of the separable quotient
problem}, dated 23 September 2026, OpenAI \cite{OpenAI} establishes two
opposite conclusions under additional set-theoretic assumptions. If the
continuum $\mathfrak c$ is real-valued measurable, every infinite-dimensional
Banach space has a separable infinite-dimensional quotient. Under the
continuum hypothesis, however, there exists a Banach space $X$ of density
$\aleph_1$ which has no separable infinite-dimensional quotient. Together
with the relative consistency of a real-valued measurable continuum from
a measurable cardinal, these results yield the independence of the
separable quotient problem from ZFC, relative to the consistency of ZFC
with a measurable cardinal.

The aim of the present paper is to describe the mechanism leading to the
space $X$ and to offer an alternative access to its definition and the
proof of its properties. We separate the functional-analytic argument
from the construction of the auxiliary space $E$, and organize the latter
in two stages. This makes it possible to see how the different properties
of $E$ enter the definition of $X$ before considering the norming sets
which produce them.

The paper is organized into four sections. Theorems~\ref{thm:abstract}
and~\ref{thm:main} are implicit in the proofs of~\cite{OpenAI}. We
formulate them explicitly here to explain how the construction yields
a space with no separable infinite-dimensional quotient.
Section~\ref{sec:goal} establishes the first criterion.

\begin{introtheorem}
Assume that a Banach space $E$ has a transfinite Schauder basis
$(e_\alpha)_{\alpha<\omega_1}$ and satisfies the following conditions.
\begin{enumerate}[label=\textup{(\roman*)},leftmargin=*]
\item Setting
\[
  X_0=\clspan\{e_\alpha^*:\alpha<\omega_1\}\subset E^*,
\]
the canonical evaluation map identifies $E$ with $X_0^*$.
\item The space $E$ contains no isomorphic copy of $\ell_1$, and every
infinite-dimensional closed subspace of $E$ has nonseparable dual.
\item There exists a closed subspace $X\subset E^*$ such that
\[
  \overline{X+X_0}^{\norm{\cdot}}=E^*
\]
and the canonical restriction map
\[
  R:E\longrightarrow X^*,\qquad (Re)(x)=x(e),
\]
is an onto isomorphism.
\end{enumerate}
Then $X$ has no infinite-dimensional separable quotient.
\end{introtheorem}

As W.~B.~Johnson pointed out, the existence of spaces satisfying all
the assumptions of this theorem implies CH. Thus CH is necessary for
this particular construction of $X$; the auxiliary space $E$ itself is
constructed in ZFC. The proof of the criterion is comparatively short.
What we find remarkable in the AI contribution is the identification of
the assumptions which make this argument possible.

Section~\ref{sec:E} gives an operator criterion for constructing $X$.

\begin{introoperator}
Let $E$ and $X_0$ satisfy conditions \textup{(i)} and \textup{(ii)}
of Theorem~\ref{thm:abstract}, with the transfinite Schauder basis
specified there. Assume that the canonical identification $E=X_0^*$
is isometric. Let $Q:E^*\to W=E^*/X_0$ be the quotient map and
$J_E:E\to E^{**}$ the canonical embedding. If a bounded linear
operator $B:E^*\to W$ satisfies
\[
  \overline{B(X_0)}^{\norm{\cdot}}=W,\qquad
  \norm{B}<\tfrac1{10},\qquad B^*(W^*)\subset J_E(E),
\]
then $X=\ker(Q-B)$ has no infinite-dimensional separable quotient.
\end{introoperator}

We then list the properties of $E$ and its special
functionals which ensure these operator hypotheses. In particular,
the adjoint range condition follows from a norm-convergence property
of $B$ on the finite initial sums of the special functionals.
Informally, $B$ may be viewed as an \emph{anti-$Q$ operator}: whereas
$Q$ annihilates $X_0$, the operator $B$ maps $X_0$ onto a norm-dense
subspace of $W$. The construction exploits these contrasting behaviours
by considering the subspace on which $Q$ and $B$ agree.

In Section~\ref{sec:construction}, we present the properties of $E$
which yield the existence of $B$. In contrast with the single construction
of the auxiliary space in \cite{OpenAI}, our description proceeds in two
stages. In the first stage we consider $E_{p,r}$, where
$1<r<p<\infty$ and $1/p+1/p'=1$. Its boundedly complete transfinite
basis $(e_\alpha)_{\alpha<\omega_1}$ admits an upper $\ell_{p'}$
estimate. Write $Q_{p,r}$ for the quotient of $E_{p,r}^*$ by the
closed span of the biorthogonal functionals. The images under $Q_{p,r}$
of the special functionals $(\phi_\gamma)_{\gamma\in\Gamma}$ admit
an upper $\ell_r$ estimate. Both estimates have constant one. Thus for
finitely supported scalar families $(a_\alpha)$ and $(b_\gamma)$,
\[
  \left\|\sum_\alpha a_\alpha e_\alpha\right\|
      \le \|(a_\alpha)\|_{\ell_{p'}},\qquad
  \left\|\sum_\gamma b_\gamma Q_{p,r}\phi_\gamma\right\|
      \le \|(b_\gamma)\|_{\ell_r}.
\]
The biorthogonal and special functionals together generate $E_{p,r}^*$
in norm, and every infinite-dimensional closed subspace of $E_{p,r}$
has nonseparable dual.

For parameters
\[
  1<r_2<p_2<r_1<p_1<\infty,
\]
put $E_i=E_{p_i,r_i}$ and $W_i=E_i^*/X_{0,i}$, with quotient maps
$Q_i:E_i^*\to W_i$. Under CH, Proposition~\ref{prop:B21} constructs
a contraction $B_{2,1}:E_2^*\to W_1$ such that
\[
  B_{2,1}e_\alpha^{(2)*}
  =Q_1\phi_{\sigma(\alpha)}^{(1)},
\]
where $\sigma$ bijectively pairs the biorthogonal functionals of $E_2$
with the indices of the special functionals of $E_1$. In particular,
$B_{2,1}(X_{0,2})$ is norm dense in $W_1$. Moreover, $B_{2,1}$ maps
the finite initial sums of every special functional of $E_2$ to a
sequence converging in norm to the image of that functional. Its
boundedness follows from the $\ell_{p_2}$ estimate for the scalar
families $(f(e_\alpha^{(2)}))_\alpha$, the upper $\ell_{r_1}$ estimate
in $W_1$, and the inclusion $\ell_{p_2}\subset\ell_{r_1}$.

In the second stage we choose
\[
  2>s>p_1>r_1>p_2>r_2>\cdots>1
\]
and describe the space
\[
  E=\left(\sum_{n=1}^{\infty}\oplus E_{p_n,r_n}\right)_{2,s}.
\]
Its dual is generated in norm by the biorthogonal functionals, the
vertical special functionals belonging to the component spaces, and
an additional family of horizontal special functionals. In $E^*/X_0$,
the vertical family at level $k$ satisfies an upper $\ell_{r_k}$
estimate and the horizontal family satisfies an upper $\ell_s$ estimate.
For $n\ge1$,
the preceding construction gives contractions
$B_{n+1,n}:E_{n+1}^*\to W_n$ because $p_{n+1}<r_n$.
The inequality $p_1<s$ similarly gives a contraction
$B_{1,0}:E_1^*\to W$ whose images of the biorthogonal functionals
generate the horizontal special cosets. Writing
$\rho_n:E^*\to E_n^*$ for restriction and $j_n:W_n\to W$ for the
natural isometric embedding, the global operator is
\[
  Bf=\tau_0 B_{1,0}(\rho_1f)
     +\sum_{n=1}^{\infty}\tau_n
          j_n B_{n+1,n}(\rho_{n+1}f),
  \qquad \tau_n=2^{-n-5}.
\]
This series converges in operator norm and gives $\norm{B}\le1/16$.
Every vertical and horizontal special coset occurs, up to a nonzero
scalar, as the image of a biorthogonal functional. Consequently
$\overline{B(X_0)}=W$. Proposition~\ref{prop:global-B} also proves
the required convergence on special functionals, and hence completes
the definition of $X=\ker(Q-B)$.

Finally, Section~\ref{sec:norming} defines the norming sets $G_{p,r}$
and $G$. An essential ingredient is the $\sigma_\varrho$ coding of
\cite[Sections~2.3--2.4]{ALT05}, based on the $\varrho$ function
of~\cite{To87}. In the separable setting an injective coding of finite
rational data is available. On $\omega_1$, this is replaced by coding
the data after transporting suitable finite $\varrho$-closures to
finite subsets of $\N$. Although $\sigma_\varrho$ is not injective,
coinciding codes give agreement on the common initial part and
$\varrho$-separation outside it. We explain how this yields the eventual
disjointness of the weight sets of two distinct infinite special
sequences. This property controls the interaction of special functionals
and is central to the definitions of the norming sets.

We prove directly several consequences of these definitions, including
bounded completeness and the estimates used in constructing $B$.
The remaining structural results, notably the assertions about
nonseparable duals of block subspaces and norm generation of the duals
by the biorthogonal and special functionals, are stated explicitly.
Their lengthy proofs use techniques from the theory of hereditarily
indecomposable Banach spaces and are outside the scope of this paper.
They may be presented separately if necessary.

Throughout the note the scalar field is $\R$. All subspaces described as
Banach spaces are norm closed, and, unless another topology is explicitly
specified, separability of a subset of a dual space refers to the norm
topology. The assertion that a space has
no separable quotient means that it has no \emph{infinite-dimensional}
separable Banach quotient.

\subsection*{The structure of \texorpdfstring{$X$}{X} and two questions}
The space $X$ has a peculiar structure. Its definition as
$\ker(Q-B)$ gives an implicit description of its vectors, and even
nonzero elements are obtained indirectly through the lifting argument.
Its subspaces and possible decompositions are therefore not transparent.
On the other hand, its dual is canonically isomorphic to the WCG space
$E$, with its explicitly defined transfinite basis and component
projections $\pi_J$. This contrast recalls the Argyros--Haydon
space~\cite{AH11}, which is hereditarily indecomposable and solves the
scalar-plus-compact problem, while its dual is isomorphic to $\ell_1$.

\begin{question}\label{question:indecomposable}
Is the space $X$ indecomposable?
\end{question}

Here indecomposability means that $X$ cannot be written as a topological
direct sum of two infinite-dimensional closed subspaces. The space $X$
has no infinite-dimensional separable complemented subspace, since a
bounded projection onto such a subspace would give a separable
infinite-dimensional quotient. Thus any decomposition into two
infinite-dimensional summands would necessarily involve two nonseparable
spaces. It is unclear whether such a decomposition is possible.
To obtain projections on $X$ from those on $E$, one needs the latter
to be weak-star continuous for the duality induced by $X$.

\begin{question}\label{question:ell-infinity}
Does there exist, under CH, a nonseparable closed subspace of
$\ell_\infty(\N)$ with no separable infinite-dimensional quotient?
\end{question}

A sufficient condition for an affirmative answer can be formulated in
terms of $E$. In what follows, the weak-star topology on $E$ is
$\sigma(E,X)$, induced by the canonical isomorphism
$R:E\to X^*$; it is the topology coming from the predual $X$.
Suppose that there is an infinite-dimensional weak-star closed subspace
$Y\subset E$ whose closed unit ball $B_Y$, with the norm inherited
from $E$, is weak-star separable. Set
\[
  M=Y_\perp=\{x\in X:x(y)=0\text{ for every }y\in Y\},
  \qquad Z=X/M.
\]
Weak-star closedness gives $M^\perp=R(Y)\subset X^*$, so that
$Z^*$ is canonically identified with $R(Y)$. In particular, $Z$ is
infinite dimensional. Choose a countable weak-star dense family
$(y_n)_n\subset B_Y$ and define
\[
  T:Z\longrightarrow\ell_\infty(\N),\qquad
  T(x+M)=(x(y_n))_{n\in\N}.
\]
This map is well defined, and weak-star density gives
\[
  \norm{T(x+M)}_\infty
  =\sup_{n\in\N}|x(y_n)|
  =\sup_{y\in B_Y}|x(y)|.
\]
Since $R:E\to X^*$ is an isomorphism, the last expression is
equivalent to the quotient norm of $x+M$. Thus $T$ is an isomorphic
embedding with closed range. Every quotient of $Z$ is also a quotient
of $X$, so $Z$ has no separable infinite-dimensional quotient.
Consequently $Z$ is nonseparable and $T(Z)$ answers
Question~\ref{question:ell-infinity} affirmatively.

There are infinite-dimensional weak-star closed subspaces $Y\subset E$
which are weak-star separable: one can take the weak-star closed linear
span of a countable linearly independent family in $E$.
However, weak-star separability of a dual space does not, in general,
imply weak-star separability of its closed unit ball. Indeed,
Avil\'es, Plebanek and Rodr\'iguez~\cite{APR14} constructed in ZFC a
compact space $K$ such that $C(K)^*$ is weak-star separable but its
closed unit ball is not. The condition on $B_Y$ therefore requires a
separate argument.

\subsection*{Acknowledgment}
This work was prepared in close collaboration with the AI model
GPT-6 Astra. Under the author's systematic supervision, the model
developed a substantial part of the statements and proofs presented
in this paper. The final manuscript was checked and edited by the
author, who takes full responsibility for its content. The author
also thanks W.~B.~Johnson for his comments and remarks.

\section{The goal}\label{sec:goal}

For a Banach space $Z$ and a subspace $Y\subset Z^*$, define the
restriction, or evaluation, operator
\[
  r_Y:Z\longrightarrow Y^*,\qquad (r_Yz)(y)=y(z).
\]
The assertion that the restriction of $Z$ to $Y$ is nonseparable means
that the set $r_Y(Z)$ is nonseparable in the norm of $Y^*$.

\begin{proposition}\label{prop:criterion}
Let $Z$ be an infinite-dimensional Banach space. Suppose that $r_Y(Z)$ is
nonseparable for every infinite-dimensional separable closed subspace
$Y\subset Z^*$. Then $Z$ has no infinite-dimensional separable quotient.
\end{proposition}

\begin{proof}
Suppose, to the contrary, that $T:Z\to S$ is a bounded linear surjection
onto an infinite-dimensional separable Banach space. Choose an
infinite-dimensional separable closed subspace $H\subset S^*$.
The adjoint $T^*:S^*\to Z^*$ is an isomorphic embedding, so
$Y=T^*(H)$ is an infinite-dimensional separable closed subspace of $Z^*$.

Define $A:S\to Y^*$ by
\[
  (As)(T^*h)=h(s)\qquad(s\in S,\ h\in H).
\]
This is well defined and bounded because $T^*|_H$ is bounded below.
For $z\in Z$ and $h\in H$,
\[
  (r_Yz)(T^*h)=(T^*h)(z)=h(Tz)=(ATz)(T^*h).
\]
Thus $r_Y=AT$ and $r_Y(Z)=A(S)$ is separable, a contradiction.
\end{proof}
\begin{theorem}\label{thm:abstract}
Assume that a Banach space $E$ has a transfinite Schauder basis
$(e_\alpha)_{\alpha<\omega_1}$ and satisfies the following conditions.
\begin{enumerate}[label=\textup{(\roman*)},leftmargin=*]
\item Setting
\[
  X_0=\clspan\{e_\alpha^*:\alpha<\omega_1\}\subset E^*,
\]
the canonical evaluation map identifies $E$ with $X_0^*$.
\item The space $E$ contains no isomorphic copy of $\ell_1$, and every
infinite-dimensional closed subspace of $E$ has nonseparable dual.
\item  There exists a closed subspace $X\subset E^*$ such that
\[
  \overline{X+X_0}^{\norm{\cdot}}=E^*
\]
and the canonical restriction map
\[
  R:E\longrightarrow X^*,\qquad (Re)(x)=x(e),
\]
is an onto isomorphism.
\end{enumerate}
Then $X$ has no infinite-dimensional separable quotient.
\end{theorem}

\begin{proof}
Let $Y\subset X^*$ be an infinite-dimensional separable closed subspace
and set $F=R^{-1}(Y)$. Since $R$ is an isomorphism, $F$ is an
infinite-dimensional separable closed subspace of $E$. Put
\[
  D_F=\{x|_F:x\in X\},\qquad
  C_F=\overline{\{y|_F:y\in X_0\}}\subset F^*.
\]

We first show that $C_F$ is separable. Every vector in $E$ is a norm
limit of finitely supported vectors, so its support is
countable. Taking the union of the supports of a countable dense subset
of $F$ gives a countable set $S\subset\omega_1$. By continuity of the
biorthogonal functionals, every vector in $F$ is supported on $S$.
Thus $e_\alpha^*|_F=0$ for $\alpha\notin S$, and
\[
  C_F=\clspan\{e_\alpha^*|_F:\alpha\in S\}
\]
is separable.

Suppose that $D_F$ were also separable. Then
$\overline{D_F+C_F}$ would be a separable closed subspace of $F^*$.
The density of $X+X_0$ in $E^*$ implies that every restriction
$f|_F$, with $f\in E^*$, belongs to $\overline{D_F+C_F}$.
By Hahn--Banach every element of $F^*$ has such an extension, whence
\[
  F^*=\overline{D_F+C_F}.
\]
This contradicts condition~\textup{(ii)}, so $D_F$ is nonseparable.

Composition with $R|_F:F\to Y$ defines an isomorphism
\[
  U:Y^*\longrightarrow F^*,\qquad U\psi=\psi\circ(R|_F).
\]
For $x\in X$ and $e\in F$,
\[
  \bigl(U(r_Yx)\bigr)(e)=(r_Yx)(Re)=(Re)(x)=x(e).
\]
Consequently $U(r_Y(X))=D_F$, and $r_Y(X)$ is nonseparable.
Proposition~\ref{prop:criterion} gives the conclusion.
\end{proof}

\begin{remark}
The isomorphism in condition~\textup{(iii)} is specified to be the
canonical restriction map: this is the compatibility used in the last
display of the proof. An abstract isomorphism $X^*\cong E$ alone does
not specify that compatibility. For $X=\ker(Q-B)$ the required
canonical isomorphism is proved in Proposition~\ref{prop:embedding}.
The argument of Theorem~\ref{thm:abstract} itself uses only the basis,
the nonseparability clause of~\textup{(ii)}, and~\textup{(iii)};
the remaining hypotheses belong to the construction.
\end{remark}

\begin{remark}
As W. B. Johnson pointed out, the existence of spaces satisfying all the
assumptions of Theorem~\ref{thm:abstract} implies CH. Here we use
density character, denoted by $\dens$, rather than algebraic dimension.
Indeed, choose an infinite-dimensional separable closed subspace $Y$ of
$E$. Its dual is nonseparable, so Stegall's theorem~\cite[Theorem~1]{S75}
gives $\dens Y^*\ge\mathfrak c$. The Hahn--Banach theorem makes the
restriction map $E^*\to Y^*$ onto, and therefore
$\dens E^*\ge\mathfrak c$.

The transfinite basis gives $\dens E=\dens X_0=\aleph_1$.
Since $X^*\cong E$, the standard inequality
$\dens X\le\dens X^*$ yields $\dens X\le\aleph_1$.
Using $\overline{X+X_0}=E^*$, we obtain
\[
  \mathfrak c\le\dens E^*
  \le\max\{\dens X,\dens X_0\}\le\aleph_1.
\]
Thus $\mathfrak c=\aleph_1$. The auxiliary space $E$ is constructed
in ZFC; CH is used to define the subspace $X$ with the
properties required in~\textup{(iii)}.
\end{remark}
\section{From some specific properties of \texorpdfstring{$E$}{E} to the space \texorpdfstring{$X$}{X}}\label{sec:E}

We first isolate the conditions on an operator $B$ that yield a space
without a separable infinite-dimensional quotient. We then give the
properties of $E$ and its special functionals which ensure these
conditions.

\subsection{An operator criterion}\label{sec:operator-criterion}

\begin{theorem}\label{thm:main}
Let $E$ be a Banach space with a transfinite Schauder basis
$(e_\alpha)_{\alpha<\omega_1}$ and biorthogonal functionals
$(e_\alpha^*)_{\alpha<\omega_1}$. Assume that:
\begin{enumerate}[label=\textup{(\roman*)},ref=\roman*,leftmargin=*]
\item\label{ass:predual}
For $X_0=\clspan\{e_\alpha^*:\alpha<\omega_1\}\subset E^*$, the
canonical evaluation map
\[
  R_0:E\longrightarrow X_0^*,\qquad (R_0e)(y)=y(e),
\]
is an onto isometry.
\item
The space $E$ contains no isomorphic copy of $\ell_1$, and every
infinite-dimensional closed subspace of $E$ has nonseparable dual.
\end{enumerate}
Let $Q:E^*\to W=E^*/X_0$ be the quotient map, where $W$ carries the
quotient norm. Suppose that there is a bounded linear operator
$B:E^*\to W$ satisfying
\begin{enumerate}[label=\textup{(B\arabic*)},ref=\textup{(B\arabic*)},leftmargin=*]
\item\label{ass:B-density}
$\overline{B(X_0)}^{\norm{\cdot}}=W$;
\item\label{ass:B-small}
$\norm{B}<1/10$;
\item\label{ass:adjoint}
$B^*(W^*)\subset J_E(E)$, where $J_E:E\to E^{**}$ is the canonical
embedding and $B^*:W^*\to E^{**}$ is the adjoint of $B$.
\end{enumerate}
Then $X=\ker(Q-B)$ has no infinite-dimensional separable quotient.
\end{theorem}

For the proof, put $L=Q-B$ and $X=\ker L$. The space $X$ is closed
in $E^*$ and carries its inherited norm. We establish three auxiliary
facts. Throughout their proofs, $\norm{B}<1/10$; the additional
assumptions are indicated where they are used.

\begin{lemma}\label{lem:lift}
Every $w\in W$ has an $L$-lift $u\in E^*$ satisfying
\[
  Lu=w,\qquad \norm{u}\le\frac52\norm{w}.
\]
\end{lemma}

\begin{proof}
For $w=0$ take $u=0$. Otherwise choose $u_0\in E^*$ with $Qu_0=w$ and
$\norm{u_0}\le2\norm{w}$. Recursively choose $u_{n+1}$ such that
\[
  Qu_{n+1}=Bu_n,\qquad
  \norm{u_{n+1}}\le2\norm{Bu_n}\le\frac15\norm{u_n}.
\]
If $Bu_n=0$, take $u_{n+1}=0$ and all subsequent terms equal to zero.
These choices are possible by the definition of the quotient norm.
The series $u=\sum_{n=0}^{\infty}u_n$ converges. Its sum satisfies
$Qu=w+Bu$, and
\[
  \norm{u}\le\frac{2}{1-1/5}\norm{w}=\frac52\norm{w}.
\]
\end{proof}

\begin{proposition}\label{prop:embedding}
Under the canonical predual and small-norm assumptions of
Theorem~\ref{thm:main}, the restriction map
\[
  R:E\longrightarrow X^*,\qquad (Re)(x)=x(e),
\]
is an isomorphic embedding, with
\[
  \frac35\norm{e}\le\norm{Re}\le\norm{e}\qquad(e\in E).
\]
If the adjoint range condition~\ref{ass:adjoint} also holds,
then $R$ is onto $X^*$.
\end{proposition}

\begin{proof}
For $y\in X_0$, apply Lemma~\ref{lem:lift} to $By$. There is $u\in E^*$
such that
\[
  Lu=By,\qquad \norm{u}\le\frac52\norm{By}\le\frac14\norm{y}.
\]
Since $Ly=-By$, the vector $x=y+u$ belongs to $X$ and satisfies
\[
  \norm{x-y}\le\frac14\norm{y},\qquad
  \norm{x}\le\frac54\norm{y}.
\]
Given $e\ne0$ and $\varepsilon>0$, the isometry $R_0$ provides
$y\in B_{X_0}$ with $|y(e)|>(1-\varepsilon)\norm{e}$. For the
corresponding $x\in X$,
\[
  \frac54\norm{Re}\ge |x(e)|
  >\left(\frac34-\varepsilon\right)\norm{e}.
\]
Let $\varepsilon$ decrease to zero. The upper bound
$\norm{Re}\le\norm{e}$ follows from the definition. Thus $R$ is
injective, has closed range, and has a bounded inverse on its range.

For surjectivity, take an arbitrary $f\in X^*$. By Hahn--Banach it has
an extension $z\in E^{**}$. The surjectivity of $R_0$ provides
$e_0\in E$ with $z|_{X_0}=J_Ee_0|_{X_0}$. Therefore $z-J_Ee_0$
annihilates $X_0$ and factors through $Q$: for some $w^*\in W^*$,
\[
  z=J_Ee_0+Q^*w^*.
\]
By condition~\ref{ass:adjoint}, choose $e_1\in E$ with $B^*w^*=J_Ee_1$.
For $x\in X$ we have $Qx=Bx$, and hence
\begin{align*}
  f(x)&=z(x)=x(e_0)+w^*(Qx)\\
      &=x(e_0)+w^*(Bx)=x(e_0+e_1).
\end{align*}
Thus $f=R(e_0+e_1)$. Since $f$ was arbitrary, $R(E)=X^*$.
\end{proof}

\begin{remark}
After $R$ has been established as an onto isomorphism, its adjoint
$R^*:X^{**}\to E^*$ is also an onto isomorphism. For $x\in X$,
\[
  R^*(J_Xx)=x\in E^*,
\]
because both sides evaluated at $e\in E$ equal $(Re)(x)=x(e)$.
This is the precise compatibility between the two dual identifications.
\end{remark}

\begin{lemma}\label{lem:density}
If $\overline{B(X_0)}=W$, then
\[
  \overline{B(X)}=W,\qquad
  \overline{Q(X)}=W,\qquad
  \overline{X+X_0}=E^*.
\]
\end{lemma}

\begin{proof}
Let $M=\overline{B(X)}$, a closed subspace of $W$. Define
\[
  A:W\longrightarrow W/M,\qquad A(Ly)=By+M\quad(y\in E^*).
\]
This is well defined: two $L$-lifts differ by an element of $X$, and
their $B$-images differ by an element of $M$.
Lemma~\ref{lem:lift} gives
\[
  \norm{A}\le\frac52\norm{B}<\frac14.
\]
For $y\in X_0$, we have $Ly=-By$, so $A(By)=-By+M$.
Since $B(X_0)$ is dense in $W$ by assumption, it follows that
$A$ is the negative of the canonical quotient map $W\to W/M$.
If $M\ne W$, that quotient map has norm one, a contradiction.
Thus $M=W$.

Since $Q=B$ on $X$, we also have $\overline{Q(X)}=W$.
Given $f\in E^*$ and $\eta>0$, choose $x\in X$ with
$\norm{Q(f-x)}<\eta$. By the quotient norm there is $y\in X_0$ with
$\norm{f-x-y}<\eta$. Hence $\overline{X+X_0}=E^*$.
\end{proof}

\begin{proof}[Proof of Theorem~\ref{thm:main}]
Proposition~\ref{prop:embedding} shows that $R:E\to X^*$ is an onto
isomorphism. Since $E$ is infinite
dimensional, so is $X$. Lemma~\ref{lem:density} gives
$\overline{X+X_0}=E^*$. Together with the assumptions on $E$, these
facts verify all the hypotheses of Theorem~\ref{thm:abstract}.
Therefore $X$ has no infinite-dimensional separable quotient.
\end{proof}

Theorem~\ref{thm:main} does not require CH as a separate hypothesis.
CH enters in the construction of an operator $B$ with the stated
properties. The theorem isolates the functional-analytic argument
from that construction.

\subsection{A list of properties of \texorpdfstring{$E$}{E}}

The following properties provide the concrete setting in which we
apply Theorem~\ref{thm:main}. Properties~\ref{ass:basis}--\ref{ass:span}
describe the space $E$ and its special functionals;
property~\ref{ass:B} describes the operator to be constructed from
this structure in Section~\ref{sec:construction}.

Let $G\subset\ell_\infty(\omega_1)$ be a norming set such that
\[
  \norm{u}_G=
  \sup_{g\in G}\left|\sum_{\alpha<\omega_1}u(\alpha)g(\alpha)\right|
  \qquad\bigl(u\in c_{00}(\omega_1)\bigr)
\]
is a finite-valued norm. Let $E$ be the completion of
$(c_{00}(\omega_1),\norm{\cdot}_G)$, and write $e_\alpha$ for its
basis vectors. The space $E$ has the following properties.

\begin{enumerate}[label=\textup{(\arabic*)},ref=\arabic*,leftmargin=*]
\item\label{ass:basis}
The family $(e_\alpha)_{\alpha<\omega_1}$ is a transfinite Schauder
basis of $E$. Its biorthogonal functionals are denoted by $e_\alpha^*$.
Setting
\[
  X_0=\clspan\{e_\alpha^*:\alpha<\omega_1\}\subset E^*,
\]
the canonical evaluation map
\[
  R_0:E\longrightarrow X_0^*,\qquad (R_0e)(y)=y(e),
\]
is an onto isometry. Thus $E=X_0^*$ through this identification.

\item\label{ass:hereditary}
The space $E$ is weakly compactly generated, contains no isomorphic copy
of $\ell_1$, and every infinite-dimensional closed subspace $F\subset E$
has nonseparable dual $F^*$.

\item\label{ass:paths}
There is a family
\[
  \Phi=\{\phi_\gamma:\gamma\in\Gamma\}\subset B_{E^*}
\]
with the following properties.
\begin{enumerate}[label=\textup{(\roman*)},leftmargin=*]
\item For each $\gamma\in\Gamma$,
\[
  \phi_\gamma=w^*\!\text{-}\!\sum_{n=1}^{\infty}\phi_{n,\gamma},
\]
where $(\phi_{n,\gamma})_n$ is a successive block sequence of finitely supported 
 functionals. The equality means convergence of the partial
sums in $\sigma(E^*,E)$.
\item If 
\[ Q: E^* \to E^*/X_0 = W \]
 is the quotient operator then \[
 Q( \phi_\gamma )\ne 0  \quad \text{ for all} \quad \gamma \in \Gamma \]

\end{enumerate}

\item\label{ass:span}
The dual space is generated in norm by the biorthogonal functionals and
the family $\Phi$:
\[
  E^*=\clspan\bigl(
    \{e_\alpha^*:\alpha<\omega_1\}\cup
    \{\phi_\gamma:\gamma\in\Gamma\}\bigr).
\]

\item\label{ass:B}
Under the continuum
hypothesis there exist a bounded linear operator
\[
  B:E^*\longrightarrow W,\qquad \norm{B}<\frac1{10},
\]
a bijection $\sigma:\omega_1\to\Gamma$, and scalars
$0<\lambda_\gamma<1$ such that
\[
  B(e_\alpha^*)=
    \lambda_{\sigma(\alpha)}Q(\phi_{\sigma(\alpha)})
    \qquad(\alpha<\omega_1).
\]
In addition, for every $\gamma\in\Gamma$,
\begin{equation}\label{eq:path-compatibility}
  B\left(\sum_{n=1}^N\phi_{n,\gamma}\right)
  =\sum_{n=1}^N B(\phi_{n,\gamma})
  \xrightarrow[N\to\infty]{\norm{\cdot}} B(\phi_\gamma).
\end{equation}
Thus the indices of the biorthogonal functionals and the path cosets are paired
bijectively. Only density of the resulting range, not surjectivity of
$B|_{X_0}$, will be used below.
\end{enumerate}

Property~\ref{ass:B} is established in
Proposition~\ref{prop:global-B} from the description of $E$ in
Section~\ref{sec:combined}.

\begin{remark}\label{rem:density}
Properties~\ref{ass:span} and~\ref{ass:B} imply
\[
  W=\clspan\{Q\phi_\gamma:\gamma\in\Gamma\}
    =\overline{B(X_0)}.
\]

Indeed, the biorthogonal functionals vanish under $Q$, and every
$Q\phi_\gamma$ is a nonzero scalar multiple of an image of a biorthogonal
functional under $B$. All closures in this formula are norm closures.
\end{remark}

\subsection{Verification of the operator hypotheses}\label{sec:duality}

Property~\ref{ass:B} gives the small-norm condition
\ref{ass:B-small}, and Remark~\ref{rem:density} gives the density
condition~\ref{ass:B-density}. It remains to verify the adjoint range
condition~\ref{ass:adjoint}. This is the role of the weak-star block
expansions and the convergence in~\eqref{eq:path-compatibility}.

\begin{proposition}\label{prop:adjoint-range}
Properties~\ref{ass:basis}, \ref{ass:paths}\textup{(i)},
\ref{ass:span}, and~\ref{ass:B} imply the adjoint range
condition~\ref{ass:adjoint} of Theorem~\ref{thm:main}:
\[
  B^*(W^*)\subset J_E(E).
\]
\end{proposition}

\begin{proof}
Fix $w^*\in W^*$. The functional
\[
  y\longmapsto w^*(By)\qquad(y\in X_0)
\]
belongs to $X_0^*$. By property~\ref{ass:basis} there is a unique $e\in E$ such
that
\begin{equation}\label{eq:represent-on-X0}
  w^*(By)=y(e)\qquad(y\in X_0).
\end{equation}
Moreover, $\norm{e}\le\norm{B}\norm{w^*}$.

For $\gamma\in\Gamma$ put
$s_{N,\gamma}=\sum_{n=1}^N\phi_{n,\gamma}$.
Each $s_{N,\gamma}$ belongs to $X_0$, since the block functionals are
finitely supported. Property~\ref{ass:paths}\textup{(i)} gives
$s_{N,\gamma}\to\phi_\gamma$ in $\sigma(E^*,E)$, whereas
\eqref{eq:path-compatibility} gives
$Bs_{N,\gamma}\to B\phi_\gamma$ in the norm of $W$. Therefore
\begin{align*}
  (B^*w^*)(\phi_\gamma)
    &=w^*(B\phi_\gamma)
      =\lim_{N\to\infty}w^*(Bs_{N,\gamma})\\
    &=\lim_{N\to\infty}s_{N,\gamma}(e)
      =\phi_\gamma(e)=(J_Ee)(\phi_\gamma).
\end{align*}
By~\eqref{eq:represent-on-X0}, the continuous functionals $B^*w^*$
and $J_Ee$ also agree on $X_0$. Property~\ref{ass:span} says that
$X_0$ and the $\phi_\gamma$ together have norm-dense linear span in
$E^*$. Hence they agree on all of $E^*$, and $B^*w^*=J_Ee$.
Since $w^*$ was arbitrary, this proves the required inclusion.
\end{proof}

\begin{remark}
Thus the convergence property~\eqref{eq:path-compatibility}
is sufficient for the adjoint range condition. The same proof would work with weak
convergence in $W$ in place of norm convergence. Neither the small
norm of $B$ nor property~\ref{ass:paths}\textup{(ii)} is needed for
this particular implication.
\end{remark}

\begin{corollary}\label{cor:main-properties}
Assume CH and properties~\textup{(1)}--\textup{(5)} above. Then
$X=\ker(Q-B)$ has no infinite-dimensional separable quotient,
$\overline{X+X_0}=E^*$, and the canonical restriction map
$R:E\to X^*$ is an onto isomorphism.
\end{corollary}

\begin{proof}
Properties~\ref{ass:basis} and~\ref{ass:hereditary} give the
assumptions on $E$ in Theorem~\ref{thm:main}.
Remark~\ref{rem:density} and property~\ref{ass:B} verify
conditions~\ref{ass:B-density} and~\ref{ass:B-small}, respectively.
Proposition~\ref{prop:adjoint-range} verifies
condition~\ref{ass:adjoint}. The conclusion follows from
Theorem~\ref{thm:main}.
\end{proof}

\section{The space \texorpdfstring{$E$}{E}}\label{sec:construction}

The auxiliary space $E$ is constructed in ZFC; CH enters
when we define the operator $B$ and the new predual $X$. The properties
needed here include density character $\aleph_1$, the absence of an
isomorphic copy of $\ell_1$, and the requirement that every
infinite-dimensional closed subspace have nonseparable dual. Thus we
seek a nonseparable analogue of the Gowers tree space $GT$~\cite{G}.

Gowers's construction combines techniques from the theory of
hereditarily indecomposable spaces with phenomena originating in the
James tree space~\cite{J}. We refer to Petsoulas and
Raikoftsalis~\cite{PR09} for a detailed study of a related Gowers-tree-like
space and its operators, and to~\cite{AAT06} for constructions using the
attractors method. The preprint~\cite{OpenAI} presents a nonseparable
space with the properties required here. I had expected that a space of
this kind could be constructed by combining the separable examples with
the methods of my joint work with J. L\'opez-Abad and
S. Todor\v{c}evi\'c~\cite{ALT05}. In particular, some years ago I
discussed a possible project in this direction with K. Beanland.

The work of Tsirel'son~\cite{T74}, Figiel and Johnson~\cite{FJ74},
Schlumprecht~\cite{S91}, and Gowers and Maurey~\cite{GM} provides
powerful methods for defining saturated norms. These methods, together
with hereditary James-tree constructions, motivate the description below.

We organize the description in two stages. First we consider spaces
$E_{p,r}$, with $1<r<p\le\infty$. In the second stage we
combine spaces $E_n=E_{p_n,r_n}$ with parameters
\[
  2>s>p_1>r_1>p_2>r_2>\cdots>1,
\]
using the norm denoted by
\[
  E=\left(\sum_{n=1}^{\infty}\oplus E_n\right)_{2,s}.
\]
This two-stage presentation differs from the single construction
in~\cite{OpenAI}. I have not checked all the proofs concerning $E$ in
that paper and therefore do not claim to have verified its construction.
I do, however, have a complete understanding of how the two spaces
are constructed.

The following is a description of the component spaces and a
proof of the operator argument that follows from their stated properties.
The norming sets defining these norms will be described in a later section.

\subsection{The spaces \texorpdfstring{$E_{p,r}$}{E(p,r)} for finite \texorpdfstring{$p$}{p}}

Fix $1<r<p<\infty$ and put $p'=p/(p-1)$. We consider the completion
of $c_{00}(\omega_1)$ under a norm $\norm{\cdot}_{p,r}$ with the
properties listed below.

\subsubsection{Properties of \texorpdfstring{$E_{p,r}$}{E(p,r)}}

\begin{enumerate}[label=\textup{(\arabic*)},leftmargin=*]
\item\label{pr:basis}
The family $(e_\alpha)_{\alpha<\omega_1}$ is a bimonotone,
boundedly complete transfinite Schauder basis of $E_{p,r}$.

\item\label{pr:wcg}
The set $K=\{0\}\cup\{e_\alpha:\alpha<\omega_1\}$ is weakly compact,
and every $e_\alpha$ is isolated in its relative weak topology.
Consequently $E_{p,r}$ is weakly compactly generated (WCG).

\item\label{pr:predual}
Writing $(e_\alpha^*)_{\alpha<\omega_1}$ for the biorthogonal
functionals and setting
\[
  X_0^{p,r}=\clspan\{e_\alpha^*:\alpha<\omega_1\}\subset E_{p,r}^*,
\]
the canonical map $E_{p,r}\to(X_0^{p,r})^*$, given by evaluation,
is an onto isometry.

\item\label{pr:hereditary}
Every infinite-dimensional closed subspace of $E_{p,r}$ has
nonseparable dual, and $E_{p,r}$ contains no isomorphic copy of $\ell_1$.

\item\label{pr:upper}
The basis vectors satisfy an upper $\ell_{p'}$ estimate with
constant one: for every finitely supported scalar family
$(a_\alpha)_{\alpha<\omega_1}$,
\begin{equation}\label{eq:pr-upper}
  \left\|\sum_{\alpha<\omega_1}a_\alpha e_\alpha\right\|_{p,r}
  \le\left(\sum_{\alpha<\omega_1}|a_\alpha|^{p'}\right)^{1/p'}.
\end{equation}
\end{enumerate}

\subsubsection{Properties of \texorpdfstring{$E_{p,r}^*$}{E(p,r)*}}

\begin{enumerate}[label=\textup{(\roman*)},leftmargin=*]
\item\label{pr:coordinates}
For every $x^*\in E_{p,r}^*$, the scalar family
$(x^*(e_\alpha))_{\alpha<\omega_1}$ belongs to $\ell_p(\omega_1)$
and satisfies
\begin{equation}\label{eq:pr-lower}
  \left(\sum_{\alpha<\omega_1}|x^*(e_\alpha)|^p\right)^{1/p}
  \le\norm{x^*}.
\end{equation}
Indeed, \eqref{eq:pr-upper} defines a contraction
\[
  T:\ell_{p'}(\omega_1)\longrightarrow E_{p,r},\qquad
  T(a)=\sum_{\alpha<\omega_1}a_\alpha e_\alpha.
\]
The series converges in norm, by approximation with finite sums.
Its adjoint is given by
\[
  T^*x^*=(x^*(e_\alpha))_{\alpha<\omega_1}\in\ell_p(\omega_1),
\]
and $\norm{T^*}\le1$, which proves \eqref{eq:pr-lower}.
In particular, the set $\{\alpha<\omega_1:x^*(e_\alpha)\ne0\}$ is countable.

This also explains property~\ref{pr:wcg}: for every $x^*\in E_{p,r}^*$
and $\varepsilon>0$, only finitely many $\alpha$ satisfy
$|x^*(e_\alpha)|\ge\varepsilon$. Hence every weak neighborhood of zero
contains all but finitely many basis vectors. It follows directly
that $K$ is weakly compact; the biorthogonal functionals isolate its
nonzero points. Since $K$ has dense linear span, $E_{p,r}$ is WCG.

\item\label{pr:paths}
There is an uncountable family of distinct functionals
\[
  \Phi_{p,r}=\{\phi_\gamma:\gamma\in\Gamma\}\subset B_{E_{p,r}^*}
\]
with the following properties.
\begin{enumerate}[label=\textup{(\alph*)},leftmargin=*]
\item\label{pr:path-expansion}
For each $\gamma\in\Gamma$ there is a successive sequence of finitely
supported blocks $(\phi_{n,\gamma})_{n=1}^{\infty}$ of the biorthogonal
functionals, with $\phi_{n,\gamma}\in B_{E_{p,r}^*}$, such that
\[
  \phi_\gamma=w^*\!\text{-}\!\sum_{n=1}^{\infty}\phi_{n,\gamma}.
\]
Thus the partial sums converge in $\sigma(E_{p,r}^*,E_{p,r})$.

\item\label{pr:path-span}
The family $\{e_\alpha^*:\alpha<\omega_1\}\cup\Phi_{p,r}$ has
norm-dense linear span in $E_{p,r}^*$.

\item\label{pr:tails}
For $\gamma\in\Gamma$ and $n\in\N$, put
\[
  \phi_{>n,\gamma}
  =\phi_\gamma-\sum_{k=1}^{n}\phi_{k,\gamma}
  =w^*\!\text{-}\!\sum_{k>n}\phi_{k,\gamma}.
\]
For every finite set $F\subset\Gamma$ there is $n_0\in\N$ such that,
for every $n>n_0$ and every scalar family $(a_\gamma)_{\gamma\in F}$,
\begin{equation}\label{eq:pr-tail-estimate}
  \left\|\sum_{\gamma\in F}a_\gamma\phi_{>n,\gamma}\right\|
  \le\left(\sum_{\gamma\in F}|a_\gamma|^r\right)^{1/r}.
\end{equation}

\item\label{pr:nonzero}
Each $\phi_\gamma$ lies outside $X_0^{p,r}$. Equivalently, if
\[
  Q_{p,r}:E_{p,r}^*\longrightarrow
  W_{p,r}:=E_{p,r}^*/X_0^{p,r}
\]
is the quotient map, then $Q_{p,r}\phi_\gamma\ne0$ for every
$\gamma\in\Gamma$. This nonvanishing is a separate requirement;
it does not follow from~\ref{pr:path-span} and~\ref{pr:tails} alone.
\end{enumerate}
\end{enumerate}

The tail estimate implies that the quotient family admits an upper
$\ell_r$ estimate with constant one. Indeed, a finite initial sum of
each path belongs to $X_0^{p,r}$, so
$Q_{p,r}\phi_{>n,\gamma}=Q_{p,r}\phi_\gamma$. For a finite
$F\subset\Gamma$, choose $n$ as in~\ref{pr:tails}. Then
\begin{equation}\label{eq:pr-quotient-estimate}
  \left\|\sum_{\gamma\in F}a_\gamma Q_{p,r}\phi_\gamma\right\|
  \le\left\|\sum_{\gamma\in F}a_\gamma\phi_{>n,\gamma}\right\|
  \le\left(\sum_{\gamma\in F}|a_\gamma|^r\right)^{1/r}.
\end{equation}
Property~\ref{pr:path-span} also gives
\[
  W_{p,r}=\clspan\{Q_{p,r}\phi_\gamma:\gamma\in\Gamma\}.
\]
Consequently the assignment of the unit vector $u_\gamma$ to
$Q_{p,r}\phi_\gamma$ extends to a contraction
\begin{equation}\label{eq:pr-synthesis}
  S_{p,r}:\ell_r(\Gamma)\longrightarrow W_{p,r},\qquad
  S_{p,r}u_\gamma=Q_{p,r}\phi_\gamma,
\end{equation}
with norm-dense range.

\subsubsection{A step towards the operator \texorpdfstring{$B$}{B}}

Section~\ref{sec:E} uses an operator that sends biorthogonal functionals
to scalar multiples of path cosets. On a single component $E_{p,r}$,
the estimates above do not by themselves establish boundedness of the
corresponding formula: the scalar family $(x^*(e_\alpha))_\alpha$ is
controlled in $\ell_p$, whereas the quotient estimate uses $\ell_r$, and $r<p$.
The formal identity from $\ell_p$ to $\ell_r$ is not bounded. This is
an obstruction to that particular estimate, not a proof that no
operator with the prescribed images of $e_\alpha^*$ can exist. The
construction resolves the exponent mismatch by mapping between
different components.

For the next definition, let
\[
  1<r_2<p_2<r_1<p_1<\infty,
\]
and abbreviate
\[
  E_i=E_{p_i,r_i},\qquad X_{0,i}=X_0^{p_i,r_i},\qquad
  W_i=E_i^*/X_{0,i}\qquad(i=1,2).
\]
Use superscripts $(i)$ on the basis vectors, their biorthogonal
functionals, and the path families to distinguish the two spaces.
Let $Q_i:E_i^*\to W_i$ denote the quotient maps.

Under CH the uncountable family $\Phi_{p_1,r_1}$ has cardinality
$\aleph_1$. To see this, note that every $x^*\in E_1^*$ is determined
by the scalar family $(x^*(e_\alpha^{(1)}))_{\alpha<\omega_1}$,
which has countable support by~\eqref{eq:pr-lower}.
Under CH there are only
$\aleph_1^{\aleph_0}=\aleph_1$ such real-valued families on $\omega_1$.
Since the path functionals are distinct and form an uncountable family,
$|\Gamma_1|=\aleph_1$.

\begin{definition}\label{def:B21}
Assume CH and fix a bijection $\sigma:\omega_1\to\Gamma_1$. For
$x^*\in E_2^*$ set
\begin{equation}\label{eq:B21-definition}
  B_{2,1}x^*
  =\sum_{\alpha<\omega_1}
    x^*(e_\alpha^{(2)})\,Q_1\phi_{\sigma(\alpha)}^{(1)}.
\end{equation}
In this particular case the sum is understood as the norm limit of the directed net of its finite subsums.
The following proposition verifies that it exists and defines a
bounded operator $B_{2,1}:E_2^*\to W_1$.
\end{definition}

\begin{proposition}\label{prop:B21}
The operator $B_{2,1}$ has the following properties.
\begin{enumerate}[label=\textup{(\arabic*)},leftmargin=*]
\item It is a well-defined bounded linear operator with
$\norm{B_{2,1}}\le1$.
\item For every $\alpha<\omega_1$,
\[
  B_{2,1}e_\alpha^{(2)*}=Q_1\phi_{\sigma(\alpha)}^{(1)},
\]
and these images have norm-dense linear span in $W_1$.
\item For every $\gamma\in\Gamma_2$,
\begin{equation}\label{eq:B21-path-compatibility}
  B_{2,1}\left(\sum_{n=1}^N\phi_{n,\gamma}^{(2)}\right)
  =\sum_{n=1}^N B_{2,1}(\phi_{n,\gamma}^{(2)})
  \xrightarrow[N\to\infty]{\norm{\cdot}}
  B_{2,1}(\phi_\gamma^{(2)}).
\end{equation}
\end{enumerate}
\end{proposition}

\begin{proof}
For $x^*\in E_2^*$ put
$a_\alpha=x^*(e_\alpha^{(2)})$. By~\eqref{eq:pr-lower},
$a=(a_\alpha)_{\alpha<\omega_1}$ belongs to $\ell_{p_2}(\omega_1)$
and $\norm{a}_{\ell_{p_2}}\le\norm{x^*}$. Since $p_2<r_1$,
\begin{equation}\label{eq:B21-coefficient-bound}
  \norm{a}_{\ell_{r_1}}
  \le\norm{a}_{\ell_{p_2}}\le\norm{x^*}.
\end{equation}
Relabel $a$ by $\sigma$ to obtain
$\widetilde a=(a_{\sigma^{-1}(\eta)})_{\eta\in\Gamma_1}$.
The contraction $S_{p_1,r_1}$ in~\eqref{eq:pr-synthesis} gives
\[
  B_{2,1}x^*=S_{p_1,r_1}\widetilde a.
\]
Finitely supported scalar families are dense in $\ell_{r_1}(\Gamma_1)$,
so this equality proves that the series
in~\eqref{eq:B21-definition} converges in norm, independently of the
order of summation. It also proves linearity and the useful estimate
\begin{equation}\label{eq:B21-bound}
  \norm{B_{2,1}x^*}
  \le\left(\sum_{\alpha<\omega_1}|x^*(e_\alpha^{(2)})|^{p_2}\right)^{1/p_2}
  \le\norm{x^*}.
\end{equation}
This proves~\textup{(1)}. Biorthogonality gives the identity
in~\textup{(2)}, and density follows because $\sigma$ is onto and
the path cosets span a dense subspace of $W_1$.

To prove~\textup{(3)}, fix $\gamma\in\Gamma_2$ and set
\[
  s_N=\sum_{n=1}^N\phi_{n,\gamma}^{(2)},\qquad
  A_N=\bigcup_{n=1}^N\supp\phi_{n,\gamma}^{(2)},\qquad
  A=\bigcup_{n=1}^{\infty}\supp\phi_{n,\gamma}^{(2)}.
\]
The block supports are pairwise disjoint. Evaluating the weak-star
expansion at each $e_\alpha^{(2)}$ therefore gives
\[
  s_N(e_\alpha^{(2)})=
  \begin{cases}
    \phi_\gamma^{(2)}(e_\alpha^{(2)}),&\alpha\in A_N,\\
    0,&\alpha\notin A_N,
  \end{cases}
  \qquad
  \phi_\gamma^{(2)}(e_\alpha^{(2)})=0\quad(\alpha\notin A).
\]
The scalar family
$(\phi_\gamma^{(2)}(e_\alpha^{(2)}))_{\alpha<\omega_1}$ belongs to
$\ell_{p_2}(\omega_1)$ by~\eqref{eq:pr-lower}, and $A_N$ increases
to $A$. Applying~\eqref{eq:B21-bound} to
$\phi_\gamma^{(2)}-s_N$, we obtain
\[
  \norm{B_{2,1}(\phi_\gamma^{(2)})-B_{2,1}s_N}
  \le\left(\sum_{\alpha\notin A_N}
       |\phi_\gamma^{(2)}(e_\alpha^{(2)})|^{p_2}\right)^{1/p_2}
  \longrightarrow0.
\]
Linearity gives the finite-sum identity
in~\eqref{eq:B21-path-compatibility}, completing the proof.
\end{proof}

\subsection{The combined space \texorpdfstring{$E$}{E}}\label{sec:combined}
In this section we present properties of the space $E =  \left(\sum_{k=1}^{\infty}\oplus E_k\right)_{2,s}$ with $E_k $ isometric to $E_{p_k,r_k}$. One of the basic properties of the space is that
every horizontal  block sequence generates a separable subspace not containg $\ell_1$ and with non separable dual. In the  separable setting this type of construction are known since 1998. In particular in 
\cite{ AF00} it is shown that for $(X_n)_{n\in \N}$ a family of separable space there exists a space $X=\left(\sum_n\oplus X_n\right)_{hi}$ with the property every horizontal block sequence generates an HI space. This type of spaces are also used in \cite{AR08}. The closest separable analogue to the space $E$ are spaces  $X=\left(\sum_n\oplus X_n\right)_{gt}$ with the property every horizontal block sequence generates an HI space with non separable dual. Such a result do not exist but it is provable.

We now proceed to  the description of $E$. Fix parameters
\begin{equation}\label{eq:combined-parameters}
  2>s>p_1>r_1>p_2>r_2>\cdots>1.
\end{equation}
In particular, $1<r_k<p_k<\infty$ for every $k$, and
$p_{k+1}<r_k$.

\begin{definition}\label{def:horizontal}
Let $C\subset\omega_1$ be unbounded and let $1<r<p<\infty$.
The set $C$ has order type $\omega_1$. We denote by $E_{p,r}(C)$
the isometric copy of $E_{p,r}$ supported on $C$, obtained by
transporting the basis along the unique order isomorphism
$C\to\omega_1$.

Fix a partition $(C_k)_{k\ge1}$ of $\omega_1$ into pairwise disjoint
unbounded sets. For an interval $J\subset\N$, put
\[
  C_J=\bigcup_{k\in J}C_k,\qquad
  \pi_J:c_{00}(\omega_1)\longrightarrow c_{00}(C_J),\qquad
  \pi_Jx=x|_{C_J},
\]
where the restriction is extended by zero outside $C_J$.
Write $\pi_k=\pi_{\{k\}}$ and, for $x\in c_{00}(\omega_1)$, put
\[
  \operatorname{ind}(x)=\{k\in\N:\pi_kx\ne0\}.
\]
A successive block sequence $(x_n)_n$ is an \emph{h-block sequence}
if the finite sets $\operatorname{ind}(x_n)$ are successive:
\[
  \max\operatorname{ind}(x_n)
  <\min\operatorname{ind}(x_{n+1})\qquad(n\in\N).
\]
For a finitely supported functional
$f=\sum_\alpha a_\alpha e_\alpha^*$, define
\[
  \operatorname{ind}(f)
  =\{k:\supp f\cap C_k\ne\varnothing\}.
\]
An h-block sequence of such functionals is defined in the same way:
the supports are successive in $\omega_1$, and their index sets are
successive in $\N$.
\end{definition}

\begin{theorem}\label{thm:combined-space}
For the partition $(C_k)_k$ and parameters
\eqref{eq:combined-parameters}, there is a norm $\norm{\cdot}_{2,s}$
on $c_{00}(\omega_1)$ such that its completion
\[
  E=\overline{\bigl(c_{00}(\omega_1),\norm{\cdot}_{2,s}\bigr)}
\]
has the following properties. Write $E_k=E_{p_k,r_k}(C_k)$.
\begin{enumerate}[label=\textup{(\arabic*)},ref=\arabic*,leftmargin=*]
\item\label{combined:basis}
The family $(e_\alpha)_{\alpha<\omega_1}$ is a boundedly complete
transfinite Schauder basis of $E$. Setting
\[
  X_0=\clspan\{e_\alpha^*:\alpha<\omega_1\},
\]
the canonical evaluation map $E\to X_0^*$ is an onto isometry.

\item\label{combined:projections}
For every interval $J\subset\N$, the map $\pi_J$ extends to a
contractive projection on $E$, with range
\[
  \pi_J(E)=\clspan\{e_\alpha:\alpha\in C_J\}.
\]
For every $k$, the range $\pi_k(E)$ is naturally isometric to $E_k$.
We use the notation
\[
  E=\left(\sum_{k=1}^{\infty}\oplus E_k\right)_{2,s}
\]
for this space and its decomposition.

\item\label{combined:hblocks}
Every normalized h-block sequence spans a subspace with nonseparable
dual, and its closed linear span contains no isomorphic copy of $\ell_1$.

\item\label{combined:hereditary}
The space $E$ contains no isomorphic copy of $\ell_1$, and every
infinite-dimensional closed subspace of $E$ has nonseparable dual.
\end{enumerate}

To describe the dual, let
\[
  \rho_k:E^*\longrightarrow E_k^*,\qquad \rho_k f=f|_{E_k},
  \qquad
  \iota_k:E_k^*\longrightarrow E^*,\qquad
  (\iota_k g)(x)=g(\pi_kx).
\]
The map $\iota_k$ is an isometry, $\rho_k$ is a contraction, and
$\pi_k^*=\iota_k\rho_k$. We identify $E_k^*$ with its image under
$\iota_k$. Then the following properties hold.
\begin{enumerate}[label=\textup{(\roman*)},leftmargin=*]
\item\label{combined:pure}
For each $k\ge1$, the family
\[
  \Phi^k=\{\phi_\gamma^k:\gamma\in\Gamma_k\}\subset B_{E_k^*}
\]
is the family from the component space $E_k$, extended to $E$ by
$\iota_k$. In particular, its members have successive finite block
expansions
\[
  \phi_\gamma^k=w^*\!\text{-}\!\sum_{n=1}^{\infty}\phi_{n,\gamma}^k,
  \qquad \supp\phi_{n,\gamma}^k\subset C_k,
\]
and
\[
  E_k^*=\clspan\bigl(\{e_\alpha^*:\alpha\in C_k\}\cup\Phi^k\bigr).
\]
Put
\[
  X_{0,k}=\clspan\{e_\alpha^*:\alpha\in C_k\},\qquad
  Q_k:E_k^*\longrightarrow W_k:=E_k^*/X_{0,k}.
\]
The nonzero cosets $Q_k\phi_\gamma^k$ have norm-dense linear span in
$W_k$ and satisfy an upper $\ell_{r_k}$ estimate with constant one.

\item\label{combined:horizontal}
There is an uncountable family of distinct functionals
\[
  \Phi^0=\{\phi_\gamma^0:\gamma\in\Gamma_0\}\subset B_{E^*}
\]
with the following properties.
\begin{enumerate}[label=\textup{(\alph*)},leftmargin=*]
\item For every $\gamma\in\Gamma_0$, there is an h-block sequence
$(\phi_{n,\gamma}^0)_{n=1}^{\infty}$ of finitely supported functionals
in $B_{E^*}$ such that
\[
  \phi_\gamma^0=w^*\!\text{-}\!\sum_{n=1}^{\infty}\phi_{n,\gamma}^0.
\]
\item Put
\[
  \phi_{>n,\gamma}^0
  =\phi_\gamma^0-\sum_{j=1}^{n}\phi_{j,\gamma}^0.
\]
For every finite $F\subset\Gamma_0$, there is $n_0\in\N$ such that,
for every $n>n_0$ and every scalar family $(a_\gamma)_{\gamma\in F}$,
\begin{equation}\label{eq:horizontal-tail}
  \left\|\sum_{\gamma\in F}a_\gamma\phi_{>n,\gamma}^0\right\|
  \le\left(\sum_{\gamma\in F}|a_\gamma|^s\right)^{1/s}.
\end{equation}
\item Each $\phi_\gamma^0$ lies outside $X_0$. Thus, for the quotient map
\[
  Q:E^*\longrightarrow W:=E^*/X_0,
\]
we have $Q\phi_\gamma^0\ne0$ for every $\gamma\in\Gamma_0$.
\end{enumerate}

\item\label{combined:span}
The dual has the norm-span description
\begin{equation}\label{eq:combined-dual-span}
  E^*=\clspan\left(\bigcup_{k\ge1}E_k^*\cup\Phi^0\right)
      =\clspan\left(\{e_\alpha^*:\alpha<\omega_1\}
                     \cup\bigcup_{k\ge0}\Phi^k\right).
\end{equation}
\end{enumerate}
\end{theorem}

We record some consequences of these properties that will be used below.
First, $P_m:=\pi_{\{1,\ldots,m\}}$ is a contraction and $P_mx\to x$
in norm for every $x\in E$. This follows first for finitely supported
vectors and then for all $x$ by density. Consequently, every $f\in E^*$
has the decomposition
\begin{equation}\label{eq:dual-weakstar-decomposition}
  f=w^*\!\text{-}\!\sum_{k=1}^{\infty}\iota_k\rho_kf,
  \qquad
  \norm{f}=\sup_m\norm{P_m^*f},\qquad
  \pi_k^*(B_{E^*})=\iota_k(B_{E_k^*}).
\end{equation}
Thus, if the notation
$E^*=(\sum_k\oplus E_k^*)_{(2,s)^*}$ is used, the component sums
are understood in this weak-star sense. Norm convergence of these sums
is not asserted; the norm-span description is
\eqref{eq:combined-dual-span}.

The space $E$ is also WCG. Indeed, each
$K_k=\{0\}\cup\{e_\alpha:\alpha\in C_k\}$ is weakly compact in
$E_k$, and hence in $E$, and $\norm{e_\alpha}\le1$. The set
\[
  K=\{0\}\cup\bigcup_{k\ge1}2^{-k}K_k
\]
is weakly compact and has dense linear span in $E$. To see compactness,
a weakly open set containing zero contains all $2^{-k}K_k$ for
sufficiently large $k$; the remaining finitely many sets are compact.

For $f\in E^*$, the component estimate~\eqref{eq:pr-lower} gives
\begin{equation}\label{eq:combined-scalar-bound}
  \left(\sum_{\alpha\in C_k}|f(e_\alpha)|^{p_k}\right)^{1/p_k}
  \le\norm{\rho_k f}\le\norm{f}\qquad(k\ge1).
\end{equation}
In particular, the set $\{\alpha:f(e_\alpha)\ne0\}$ is a countable
union of countable sets and is therefore countable.

\subsection{The structure of the quotient \texorpdfstring{$W$}{W}}\label{sec:quotient}

We first place all the path cosets in the common quotient $W$.
The projection $\pi_k^*$ maps $X_0$ onto $X_{0,k}$, so
\[
  E_k^*\cap X_0=X_{0,k}.
\]
Consequently the map
\[
  j_k:W_k\longrightarrow W,\qquad j_k(Q_kg)=Q(\iota_kg),
\]
is an isometric embedding. Indeed, extension by $\iota_k$ gives one
inequality between the quotient norms, and restriction by $\rho_k$
gives the reverse inequality. In particular, every $Q\phi_\gamma^k$
is nonzero, and the upper $\ell_{r_k}$ estimate also holds in $W$.

We identify each $W_k$, $k\ge1$, with its image under $j_k$, and set
\[
  W_0=\clspan\{Q\phi_\gamma^0:\gamma\in\Gamma_0\},
  \qquad
  V=\clspan\left(\bigcup_{k\ge1}W_k\right).
\]
Thus $W_0$ is the closed span of the horizontal special cosets, and
$V$ is the closed span of the vertical quotient spaces.

\begin{proposition}[The decomposition of the quotient]
\label{prop:quotient-decomposition}
The quotient $W=E^*/X_0$ has the following properties.
\begin{enumerate}[label=\textup{(\arabic*)},leftmargin=*]
\item For every $k\ge1$, the formula
\[
  \mathcal P_k(Qf)=Q(\pi_k^*f)\qquad(f\in E^*)
\]
defines a contractive projection from $W$ onto $W_k$.
\item We have $W=W_0\oplus V$, and both associated projections
$\mathcal P_0:W\to W_0$ and $\mathcal P_V:W\to V$ are
contractive. More precisely,
\begin{equation}\label{eq:quotient-separation}
  \norm{v-z}\ge\max\{\norm{v},\norm{z}\}
  \qquad(v\in W_0,\ z\in V).
\end{equation}
In particular, any two normalized vectors, one in $W_0$ and the other
in $V$, have distance at least one.
\item The family $(W_k)_{k\ge0}$ is a Schauder decomposition of $W$:
every $w\in W$ has a unique representation
\begin{equation}\label{eq:quotient-Schauder-decomposition}
  w=w_0+\sum_{k=1}^{\infty}w_k,
  \qquad w_k=\mathcal P_kw\in W_k\quad(k\ge0),
\end{equation}
where the series converges in norm. Thus we may write
\[
  W=\left(\sum_{k=0}^{\infty}\oplus W_k\right)
\]
with its inherited quotient norm. The partial-sum projections
$\mathcal P_0+\sum_{k=1}^m\mathcal P_k$ have norm at most two.
\end{enumerate}
\end{proposition}

\begin{proof}
Recall that $P_m=\pi_{\{1,\ldots,m\}}$. The operators $P_m$ and
$T_m=I-P_m=\pi_{\{m+1,m+2,\ldots\}}$ are contractions on $E$.
Their adjoints preserve $X_0$, since they send each biorthogonal
functional either to itself or to zero. They therefore induce
contractions on the quotient:
\[
  \widehat P_m(Qf)=Q(P_m^*f),\qquad
  \widehat T_m(Qf)=Q(T_m^*f).
\]
In particular, $\widehat T_m=I_W-\widehat P_m$.
The same argument applied to $\pi_k$ gives the contractive projection
$\mathcal P_k$ in~\textup{(1)}. Its range is $W_k$, and
\[
  \widehat P_m=\sum_{k=1}^m\mathcal P_k.
\]

For every horizontal special functional $\phi_\gamma^0$, we have
\begin{equation}\label{eq:horizontal-finite-component-restriction}
  P_m^*\phi_\gamma^0\in X_0.
\end{equation}
Indeed, the expansion
$\phi_\gamma^0=w^*\!\text{-}\!\sum_n\phi_{n,\gamma}^0$
is an h-block expansion, so only finitely many of its blocks meet the
first $m$ components. By weak-star continuity of $P_m^*$, for some $N$,
\[
  P_m^*\phi_\gamma^0
  =\sum_{n=1}^{N}P_m^*\phi_{n,\gamma}^0.
\]
The right-hand side is finitely supported and hence belongs to $X_0$.
Passing to the quotient and then to the closed linear span gives
\begin{equation}\label{eq:horizontal-quotient-projections}
  \widehat P_m v=0,\qquad \widehat T_m v=v
  \qquad(v\in W_0).
\end{equation}

If $z\in W_1+\cdots+W_m$, then
$\widehat P_m z=z$ and $\widehat T_m z=0$.
Consequently, for $v\in W_0$,
\[
  \norm{v-z}\ge\norm{\widehat T_m(v-z)}=\norm{v},
  \qquad
  \norm{v-z}\ge\norm{\widehat P_m(v-z)}=\norm{z}.
\]
The union of the finite vertical sums is dense in $V$, so these
inequalities extend to every $z\in V$ and give
\eqref{eq:quotient-separation}.

Replacing $z$ by $-z$, we see that $W_0\cap V=\{0\}$ and that the
two natural projections on $W_0+V$ are contractive. The sum is closed:
if $(v_n+z_n)_n$ is Cauchy, with $v_n\in W_0$ and $z_n\in V$,
then both $(v_n)_n$ and $(z_n)_n$ are Cauchy by the same inequality,
and their limits belong to the respective closed subspaces.
The dual norm-span description~\eqref{eq:combined-dual-span} implies
that $W_0+V$ is dense in $W$. Thus $W=W_0\oplus V$, proving
\textup{(2)}.

For every $z\in V$ we have $\widehat P_m z\to z$ in norm.
This holds eventually with equality on each finite vertical sum and
extends to $V$ by density and the uniform bound
$\norm{\widehat P_m}\le1$. Write $w=v+z$, with $v\in W_0$ and
$z\in V$. Equation~\eqref{eq:horizontal-quotient-projections} gives
\[
  \mathcal P_0w+\sum_{k=1}^m\mathcal P_kw
  =v+\widehat P_m z\longrightarrow v+z=w.
\]
Each $\mathcal P_k$ is the identity on $W_k$ and vanishes on the
other summands: for $k\ge1$ this follows from the component
projections and~\eqref{eq:horizontal-quotient-projections}, and
for $k=0$ it follows from~\textup{(2)}. Applying these bounded
projections to any norm-convergent representation proves uniqueness.
Finally,
\[
  \left\|\mathcal P_0+\sum_{k=1}^m\mathcal P_k\right\|
  =\norm{\mathcal P_0+\widehat P_m}\le2,
\]
which completes the proof of~\textup{(3)}. The projections onto the
horizontal and vertical parts can also be expressed by the norm limits
\[
  \mathcal P_0w=\lim_{m\to\infty}\widehat T_mw,
  \qquad
  \mathcal P_Vw=\lim_{m\to\infty}\widehat P_mw
  \qquad(w\in W).
\]
\end{proof}

This argument uses the contractive component projections, the h-block
expansions of the horizontal special functionals, and the dual norm-span
description. It does not use CH.

\subsection{The operator \texorpdfstring{$B$}{B}}\label{sec:operator-B}

For the horizontal family, finite initial block sums lie in $X_0$.
Hence \eqref{eq:horizontal-tail}, followed by $Q$, gives an upper
$\ell_s$ estimate for $(Q\phi_\gamma^0)_{\gamma\in\Gamma_0}$.
Set
\[
  q_0=s,\qquad q_k=r_k\quad(k\ge1).
\]
The argument proving~\eqref{eq:pr-synthesis} now gives contractions
\begin{equation}\label{eq:global-synthesis}
  S_k:\ell_{q_k}(\Gamma_k)\longrightarrow W,\qquad
  S_ku_\gamma=Q\phi_\gamma^k\qquad(k\ge0),
\end{equation}
where $(u_\gamma)_\gamma$ denotes the unit vector basis. Also,
\eqref{eq:combined-dual-span} implies
\begin{equation}\label{eq:all-path-cosets}
  W=\clspan\{Q\phi_\gamma^k:k\ge0,\ \gamma\in\Gamma_k\}.
\end{equation}

Assume CH. Every functional in $E^*$ is determined by a countably
supported real-valued family on $\omega_1$, by
\eqref{eq:combined-scalar-bound}. There are
$\aleph_1^{\aleph_0}=\aleph_1$ such families under CH.
Each $\Phi^k$ is uncountable and consists of distinct functionals;
therefore $|\Gamma_k|=\aleph_1$ for every $k\ge0$. We may fix bijections
\begin{equation}\label{eq:global-bijections}
  \sigma_k:C_{k+1}\longrightarrow\Gamma_k\qquad(k\ge0).
\end{equation}
Thus $C_1$ is assigned to the horizontal family $\Phi^0$, and
$C_{k+1}$ is assigned to $\Phi^k$ for $k\ge1$. The parameter inequalities
give exactly the estimates needed for these assignments:
\[
  p_1<s=q_0,\qquad p_{k+1}<r_k=q_k\quad(k\ge1).
\]

Figure~\ref{fig:operator-B-families} illustrates the assignment of the
families and the component operators in the following proposition.

\begin{proposition}\label{prop:global-B}
Assume CH and the properties of $E$ in
Theorem~\ref{thm:combined-space}. For $k\ge0$ put $\tau_k=2^{-k-5}$.
For $f\in E^*$ define
\begin{equation}\label{eq:global-B}
  Bf=\sum_{k=0}^{\infty}\tau_k S_k
       \bigl((f(e_{\sigma_k^{-1}(\gamma)}))_{\gamma\in\Gamma_k}\bigr).
\end{equation}
Then $B:E^*\to W$ is a bounded linear operator satisfying:
\begin{enumerate}[label=\textup{(\arabic*)},leftmargin=*]
\item $\norm{B}\le1/16<1/10$.
\item For $k\ge0$ and $\alpha\in C_{k+1}$,
\begin{equation}\label{eq:global-B-basis}
  B(e_\alpha^*)=\tau_k Q\phi_{\sigma_k(\alpha)}^k.
\end{equation}
In particular, $\overline{B(X_0)}=W$.
\item For every $j\ge0$ and $\gamma\in\Gamma_j$,
\begin{equation}\label{eq:global-B-convergence}
  B\left(\sum_{n=1}^{N}\phi_{n,\gamma}^j\right)
  =\sum_{n=1}^{N}B(\phi_{n,\gamma}^j)
  \xrightarrow[N\to\infty]{\norm{\cdot}}B(\phi_\gamma^j).
\end{equation}
\item $B^*(W^*)\subset J_E(E)$.
\end{enumerate}
Moreover, with the disjoint index set
\[
  \Gamma=\{(k,\gamma):k\ge0,\ \gamma\in\Gamma_k\},
  \qquad \phi_{(k,\gamma)}=\phi_\gamma^k,
\]
the rule $\sigma(\alpha)=(k,\sigma_k(\alpha))$ for
$\alpha\in C_{k+1}$ defines a bijection $\sigma:\omega_1\to\Gamma$.
Taking $\lambda_{(k,\gamma)}=\tau_k$, the operator $B$ satisfies
property~\ref{ass:B} of Section~\ref{sec:E}.
\end{proposition}

\clearpage
\begin{figure}[p]
\centering
\begingroup
\setlength{\unitlength}{1mm}
\small

\noindent\textbf{(a) Assignment of the families}\par
\smallskip
\begin{picture}(152,40)
  \put(13,34){\makebox(0,0){$\{e_\alpha^*:\alpha\in C_1\}$}}
  \put(40,34){\makebox(0,0){$\{e_\alpha^*:\alpha\in C_2\}$}}
  \put(67,34){\makebox(0,0){$\{e_\alpha^*:\alpha\in C_3\}$}}
  \put(85,34){\makebox(0,0){$\cdots$}}
  \put(104,34){\makebox(0,0){$\{e_\alpha^*:\alpha\in C_k\}$}}
  \put(132,34){\makebox(0,0){$\{e_\alpha^*:\alpha\in C_{k+1}\}$}}
  \put(150,34){\makebox(0,0){$\cdots$}}

  \put(13,28){\vector(0,-1){15}}
  \put(40,28){\vector(0,-1){15}}
  \put(67,28){\vector(0,-1){15}}
  \put(104,28){\vector(0,-1){15}}
  \put(132,28){\vector(0,-1){15}}
  \put(15,21){\makebox(0,0)[l]{$\scriptstyle\sigma_0$}}
  \put(42,21){\makebox(0,0)[l]{$\scriptstyle\sigma_1$}}
  \put(69,21){\makebox(0,0)[l]{$\scriptstyle\sigma_2$}}
  \put(106,21){\makebox(0,0)[l]{$\scriptstyle\sigma_{k-1}$}}
  \put(134,21){\makebox(0,0)[l]{$\scriptstyle\sigma_k$}}

  \put(13,7){\makebox(0,0){$\{\phi_\gamma^0:\gamma\in\Gamma_0\}$}}
  \put(40,7){\makebox(0,0){$\{\phi_\gamma^1:\gamma\in\Gamma_1\}$}}
  \put(67,7){\makebox(0,0){$\{\phi_\gamma^2:\gamma\in\Gamma_2\}$}}
  \put(85,7){\makebox(0,0){$\cdots$}}
  \put(104,7){\makebox(0,0){$\{\phi_\gamma^{k-1}:\gamma\in\Gamma_{k-1}\}$}}
  \put(132,7){\makebox(0,0){$\{\phi_\gamma^k:\gamma\in\Gamma_k\}$}}
  \put(150,7){\makebox(0,0){$\cdots$}}
\end{picture}

Under $CH$ the bijection $\sigma_k:C_{k+1}\to\Gamma_k$  assigns
$\phi_{\sigma_k(\alpha)}^k$ to $e_\alpha^*$. The actual image under $B$ is
\[
  B(e_\alpha^*)=\tau_kQ\phi_{\sigma_k(\alpha)}^k
  \qquad(\alpha\in C_{k+1},\ k\ge0).
\]

\medskip
\noindent\textbf{(b) The component spaces and operators}\par
\smallskip
\begin{picture}(152,37)
  \put(13,32){\makebox(0,0){$E_1^*$}}
  \put(40,32){\makebox(0,0){$E_2^*$}}
  \put(67,32){\makebox(0,0){$E_3^*$}}
  \put(85,32){\makebox(0,0){$\cdots$}}
  \put(104,32){\makebox(0,0){$E_k^*$}}
  \put(132,32){\makebox(0,0){$E_{k+1}^*$}}
  \put(150,32){\makebox(0,0){$\cdots$}}

  \put(13,26){\vector(0,-1){14}}
  \put(40,26){\vector(0,-1){14}}
  \put(67,26){\vector(0,-1){14}}
  \put(104,26){\vector(0,-1){14}}
  \put(132,26){\vector(0,-1){14}}
  \put(15,19){\makebox(0,0)[l]{$\scriptstyle B_{1,0}$}}
  \put(42,19){\makebox(0,0)[l]{$\scriptstyle B_{2,1}$}}
  \put(69,19){\makebox(0,0)[l]{$\scriptstyle B_{3,2}$}}
  \put(106,19){\makebox(0,0)[l]{$\scriptstyle B_{k,k-1}$}}
  \put(134,19){\makebox(0,0)[l]{$\scriptstyle B_{k+1,k}$}}

  \put(13,6){\makebox(0,0){$W_0$}}
  \put(40,6){\makebox(0,0){$W_1$}}
  \put(67,6){\makebox(0,0){$W_2$}}
  \put(85,6){\makebox(0,0){$\cdots$}}
  \put(104,6){\makebox(0,0){$W_{k-1}$}}
  \put(132,6){\makebox(0,0){$W_k$}}
  \put(150,6){\makebox(0,0){$\cdots$}}
\end{picture}

Here $E_k=E_{p_k,r_k}(C_k)$ and
\[
  W_0=\overline{\operatorname{span}}\{Q\phi_\gamma^0:\gamma\in\Gamma_0\}
  \subset W=E^*/X_0,
  \qquad W_k=E_k^*/X_{0,k}\quad(k\ge1).
\]
For $k\ge1$, $B_{k+1,k}$ is the operator of Proposition~\ref{prop:B21};
$B_{1,0}$ is its horizontal analogue, defined by
\[
  B_{1,0}g=\sum_{\alpha\in C_1}g(e_\alpha)\,
      Q\phi_{\sigma_0(\alpha)}^0\qquad(g\in E_1^*).
\]
Let $\rho_k:E^*\to E_k^*$ be restriction, let $j_k:W_k\to W$ be the
canonical embedding for $k\ge1$, and let $j_0:W_0\hookrightarrow W$
be inclusion. The component maps combine as
\[
  \boxed{\displaystyle
  B=\sum_{k=0}^{\infty}\tau_k\,
       \underbrace{j_kB_{k+1,k}\rho_{k+1}}_{A_k},
  \qquad \tau_k=2^{-k-5}.}
\]
\endgroup
\caption{The assignment and the operators in Proposition~\ref{prop:global-B}.
The first component supplies the horizontal family, and component
$k+1$ supplies family $k$. The weights are applied when the component
operators are combined to form $B$.}
\label{fig:operator-B-families}
\end{figure}
\clearpage

\begin{proof}
For $k\ge0$ and $f\in E^*$, write
\[
  a_k(f)=(f(e_{\sigma_k^{-1}(\gamma)}))_{\gamma\in\Gamma_k},
  \qquad A_kf=S_ka_k(f).
\]
By~\eqref{eq:combined-scalar-bound} and $p_{k+1}<q_k$,
\begin{equation}\label{eq:global-A-bound}
  \norm{A_kf}
  \le\norm{a_k(f)}_{\ell_{q_k}}
  \le\norm{a_k(f)}_{\ell_{p_{k+1}}}
  =\left(\sum_{\alpha\in C_{k+1}}|f(e_\alpha)|^{p_{k+1}}\right)^{1/p_{k+1}}
  \le\norm{f}.
\end{equation}
This is precisely the estimate of Proposition~\ref{prop:B21}.
For $k\ge1$, $A_k$ is the composition of restriction to $E_{k+1}$,
the operator from that proposition with target $W_k$, and the embedding
$j_k$. For $k=0$, the same proof uses the upper $\ell_s$ estimate of
the horizontal family in place of the upper $\ell_{r_k}$ estimate.

Each $A_k$ is therefore a linear contraction. Since
\[
  \sum_{k=0}^{\infty}\tau_k=\frac1{16},
\]
the series $B=\sum_{k=0}^{\infty}\tau_kA_k$ converges in operator
norm and proves~\textup{(1)}. Biorthogonality and the disjointness of
the $C_k$ give~\eqref{eq:global-B-basis}. Every path coset is thus a
nonzero scalar multiple of $B(e_\alpha^*)$ for a suitable $\alpha$.
Equation~\eqref{eq:all-path-cosets} yields
$\overline{B(X_0)}=W$, proving~\textup{(2)}. The global bijection and
the scalars $\lambda_{(k,\gamma)}$ have the asserted properties by
their definitions.

We prove~\textup{(3)} simultaneously for the component and horizontal
families. Fix $j\ge0$ and $\gamma\in\Gamma_j$, and put
\[
  \phi=\phi_\gamma^j,\qquad
  s_N=\sum_{n=1}^{N}\phi_{n,\gamma}^j,\qquad
  F_N=\bigcup_{n=1}^{N}\supp\phi_{n,\gamma}^j.
\]
The block supports are pairwise disjoint. The weak-star expansion
therefore gives
\[
  (\phi-s_N)(e_\alpha)=
  \begin{cases}
    0,&\alpha\in F_N,\\
    \phi(e_\alpha),&\alpha\notin F_N,
  \end{cases}
  \qquad
  \phi(e_\alpha)=0\quad
  \left(\alpha\notin\bigcup_NF_N\right).
\]
For each fixed $k$, the scalar family
$(\phi(e_\alpha))_{\alpha\in C_{k+1}}$ belongs to
$\ell_{p_{k+1}}(C_{k+1})$. Hence~\eqref{eq:global-A-bound} gives
\begin{equation}\label{eq:global-path-tail}
  \norm{A_k(\phi-s_N)}
  \le\left(\sum_{\alpha\in C_{k+1}\setminus F_N}
            |\phi(e_\alpha)|^{p_{k+1}}\right)^{1/p_{k+1}}
  \longrightarrow0.
\end{equation}
This is the same truncation argument as in
Proposition~\ref{prop:B21}.

Since $s_N\to\phi$ weak-star, the uniform boundedness principle gives
$M:=\sup_N\norm{\phi-s_N}<\infty$. For every $K\ge0$,
\[
  \norm{B(\phi-s_N)}
  \le\sum_{k=0}^{K}\tau_k\norm{A_k(\phi-s_N)}
     +M\sum_{k>K}\tau_k.
\]
For fixed $K$, the finite sum tends to zero by
\eqref{eq:global-path-tail}; the second term tends to zero as
$K\to\infty$. This proves the norm convergence in
\eqref{eq:global-B-convergence}. The finite-sum identity follows from
linearity.

Finally, Theorem~\ref{thm:combined-space} supplies the canonical
identification $E=X_0^*$, the weak-star block expansions, and the
norm-span assertion~\eqref{eq:combined-dual-span}. Together with
\eqref{eq:global-B-convergence}, these are exactly the hypotheses of
Proposition~\ref{prop:adjoint-range}. It follows that
$B^*(W^*)\subset J_E(E)$, proving~\textup{(4)}.
\end{proof}

The description of $E$, the WCG argument above, and
Proposition~\ref{prop:global-B} establish all the properties used in
Section~\ref{sec:E}. Therefore Theorem~\ref{thm:main} applies: under
CH, $X=\ker(Q-B)$ has no infinite-dimensional separable quotient,
and the canonical restriction map $E\to X^*$ is an onto isomorphism.

\section{The norming sets $G_{p,r}$ and $G$ of $E_{p,r}$ and $E$}\label{sec:norming}
In this section we will describe the norming sets for defining the norms of
the spaces $E_{p,r}$. and $E$. A norming set is a $G\subset (c_{00}(\omega_1), \Q )$
containing the singletons and for every $f\in G$ we have $\|f\|_\infty \le 1$.
It induces a norm on $ c_{00}(\omega_1) $ with the rule
\[ \|x\|_G = sup\{ |f(x)| : f\in G\} \]
If the set $G$ is closed in the projection  of its elements to the intervals 
of $\omega_1$ then the family $\{e_\alpha\}_{\alpha\in \omega_1}$ is a transfinite
basis for the space $X_G = {\clspan}  c_{00}(\omega_1, \| \cdot\|_G) $ .

For the norming sets we fix two strictly increasing sequences of
positive integers $(m_j)_j$ and $(n_j)_j$, following the
rapid-growth parameter choices used in~\cite{DM07}. We record
explicitly the conditions needed in the arguments below:
\begin{equation}\label{eq:parameter-growth}
 \sum_{j=1}^{\infty}\frac1{m_j}<1,
 \qquad
 \frac{n_j^{1/q}}{m_j}\longrightarrow\infty
 \quad\text{for every }1<q<\infty.
\end{equation}
The second condition strengthens $n_j/m_j\to\infty$ and can be
arranged, for example, by choosing $n_j\ge m_j^j$.
For bounded completeness at a fixed $p$, we use it with
$q=p'=p/(p-1)$. The standard parameters in
\cite[Section~2.1]{ALT05} also satisfy this growth condition.
\begin{definition}\label{def:pj-operations}
Let $G\subset c_{00}(\omega_1,\Q)$, $1<p<\infty$, and $j\in\N$.
\begin{enumerate}
\item We say that an $f\in G$ is the result of the $(p, j)$- operation if 
\[ f = \frac{1}{m_j}\sum_{k=1}^{n_j }\lambda_k f_k \quad  \text{with} \quad  (\lambda_k)_k \subset \Q \quad  \sum_k |\lambda_k|^p \le 1\, \text{and }\,  f_1 < \cdots < f_{n_j} \]
where $f_1,\ldots,f_{n_j}$ are nonzero elements of $G$.
\item If $ f\in G$ is the result of $(p,j)$-operation then we say that $f$ is a weighted  functional with its weight $w(f) = m_j$. Notice that for a weighted functional $f$ its weight is not necessarily unique.
\item The set $G \subset (c_{00}(\omega_1), \Q )$ is closed in $(p, j)$-operation if
\[ \forall(f_k)_{k=1}^{n_j}  \quad \text{successive blocks in} \; G \quad \forall(\lambda_k) \subset \Q \quad \text{with} \; \sum_k |\lambda_k |^p \le 1\quad \text{we have}\]
\[ f= \frac{1}{m_j} \sum_{k=1}^{n_j} \lambda_k f_k \in G\]
\end{enumerate}

\end{definition}

\subsection{The \texorpdfstring{$\sigma_\varrho$}{sigma-rho} coding and special sequences}
Special sequences are defined using the $\sigma_\varrho$ coding of
\cite[Sections~2.3--2.4]{ALT05}, based on a $\varrho$ function
as in~\cite{To87}. The coding assigns a weight $m_j$ to a finite
sequence of successive rational functionals together with their chosen
weights and integer parameters. In the separable setting the family
to be coded is countable, so an injective coding is available.
On $\omega_1$, the function $\sigma_\varrho$ is not injective;
the finite closures associated with $\varrho$ provide the agreement
on common initial parts that will be used below.

\subsubsection{The \texorpdfstring{$\sigma_\varrho$}{sigma-rho} coding}
We give the construction in the following five steps. The auxiliary
function $\sigma$ takes integer values, whereas $\sigma_\varrho$ will
take values in $\{m_j:j\in\N\}$, in accordance with our convention
for weights. Thus our $\sigma_\varrho$ includes the passage from an
index to the corresponding weight; in~\cite[Section~2.4]{ALT05}
this passage is written separately.

\paragraph{1. The $\varrho$ function.}
Fix a function $\varrho:[\omega_1]^2\to\N$ with the following
properties:
\begin{enumerate}[label=\textup{(\roman*)}]
\item For $\alpha<\beta<\gamma<\omega_1$,
\[
 \varrho(\alpha,\gamma)
 \le\max\{\varrho(\alpha,\beta),\varrho(\beta,\gamma)\}.
\]
\item For $\alpha<\beta<\gamma<\omega_1$,
\[
 \varrho(\alpha,\beta)
 \le\max\{\varrho(\alpha,\gamma),\varrho(\beta,\gamma)\}.
\]
\item For every $\beta<\omega_1$ and $p\in\N$, the set
\[
 \{\alpha<\beta:\varrho(\alpha,\beta)\le p\}
\]
is finite.
\end{enumerate}
We also put $\varrho(\alpha,\alpha)=0$. The existence of such a
function is established in~\cite{To87}; we fix one throughout.

\paragraph{2. Finite closures.}
For a finite set $F\subset\omega_1$, put
\[
 p_F=\max\bigl(\{\varrho(\alpha,\beta):
                 \alpha,\beta\in F,\ \alpha<\beta\}\cup\{0\}\bigr).
\]
For $p\in\N$, define its \emph{$p$-closure} by
\[
 \overline F^{\,p}
 =F\cup\{\alpha<\omega_1:
           \text{there is }\beta\in F\text{ with }
           \alpha<\beta\text{ and }\varrho(\alpha,\beta)\le p\}.
\]
This set is finite by~\textup{(iii)}. We call $F$ \emph{$p$-closed}
if $\overline F^{\,p}=F$, and \emph{$\varrho$-closed} if
$\overline F^{\,p_F}=F$, using the same formula when $p_F=0$.
The closure is monotone and idempotent:
\[
 F\subset H\ \Longrightarrow\
 \overline F^{\,p}\subset\overline H^{\,p},
 \qquad
 \overline{\bigl(\overline F^{\,p}\bigr)}^{\,p}
 =\overline F^{\,p}.
\]
Indeed, if $\alpha\le\beta\le\gamma$, with
$\varrho(\alpha,\beta)\le p$ and $\varrho(\beta,\gamma)\le p$,
then~\textup{(i)} gives $\varrho(\alpha,\gamma)\le p$.

If $F\ne\varnothing$ and $p\ge p_F$, then
\begin{equation}\label{eq:rho-singleton-closure}
 \overline F^{\,p}=\overline{\{\max F\}}^{\,p},
 \qquad p_{\overline F^{\,p}}\le p.
\end{equation}
For the equality, use $F\subset\overline{\{\max F\}}^{\,p}$
and idempotence. For the inequality, apply~\textup{(ii)} to two
members of this closure and $\max F$.

Here is the coherence property needed for coding. If $H$ and $K$
are finite $p$-closed sets with $p_H,p_K\le p$, then for every
$\alpha\in H\cap K$,
\begin{equation}\label{eq:rho-closure-coherence}
 H\cap[0,\alpha]
 =\overline{\{\alpha\}}^{\,p}
 =K\cap[0,\alpha].
\end{equation}
In fact, $p_H\le p$ gives
$H\cap[0,\alpha]\subset\overline{\{\alpha\}}^{\,p}$;
the reverse inclusion follows from $p$-closedness. The argument
for $K$ is the same. Consequently, $H\cap K$ is an initial part
of both $H$ and $K$ (possibly empty).

\paragraph{3. Increasing enumerations of finite sets.}
For a nonempty finite set $F\subset\omega_1$, let
\[
 \pi_F:\{1,\ldots,|F|\}\longrightarrow F
\]
be its increasing bijection. If $f$ is a scalar-valued function
on $F$, write $\pi_F(f)$ for the function on $\N$ given by
\[
 \pi_F(f)(n)=
 \begin{cases}
 f(\pi_F(n)),&1\le n\le |F|,\\
 0,&n>|F|.
 \end{cases}
\]
Thus $\pi_F(f)$ records the values of $f$ in the order of $F$,
including the zero entries. We use the same notation if $f$ is
a finitely supported function on $\omega_1$ with $\supp f\subset F$.
In particular, rational-valued $f$ gives
$\pi_F(f)\in c_{00}(\N,\Q)$.

\paragraph{4. The countable coding $\sigma$.}
Let $\mathcal Q_s$ be the set of all finite sequences
\[
 s=(g_1,w_1,\ell_1,\ldots,g_d,w_d,\ell_d),\qquad d\ge1,
\]
where $g_1<\cdots<g_d$ are nonzero elements of $c_{00}(\N,\Q)$,
and $(w_i)_{i=1}^d$ and $(\ell_i)_{i=1}^d$ are strictly increasing
sequences of positive integers. Put
\[
 \varepsilon(s)=\min\{|g_i(n)|:
                    1\le i\le d,\ n\in\supp g_i\}>0.
\]
Since $\mathcal Q_s$ is countable, we may choose an injection
$\sigma:\mathcal Q_s\to\N$ satisfying
\begin{equation}\label{eq:sigma-growth}
 \sigma(s)>
 \max\{w_d,\ell_d^2,\varepsilon(s)^{-2},\max\supp g_d\}.
\end{equation}
To do so, enumerate $\mathcal Q_s$ and assign to each successive
element an unused integer exceeding its displayed bound.
The additional bound by $w_d$ makes the increase of the resulting
weights explicit.

\paragraph{5. The coding $\sigma_\varrho$.}
Let $\mathcal Q_s(\omega_1)$ consist of the finite sequences
\[
 \mathbf f=(f_1,w_1,\ell_1,\ldots,f_d,w_d,\ell_d),
\]
where $f_1<\cdots<f_d$ are nonzero elements of
$c_{00}(\omega_1,\Q)$, the sequences $(w_i)$ and $(\ell_i)$ are
strictly increasing sequences of positive integers, and
\[
 \ell_i\ge p_{F_i},\qquad
 F_i=\bigcup_{k=1}^{i}\supp f_k\quad(1\le i\le d).
\]
The parameters $\ell_i$ are chosen subject to these conditions;
they are not uniquely determined by the supports.
Set
\[
 H_{\mathbf f}=\overline{F_d}^{\,\ell_d},\qquad
 \widehat{\mathbf f}
 =\bigl(\pi_{H_{\mathbf f}}(f_1),w_1,\ell_1,\ldots,
         \pi_{H_{\mathbf f}}(f_d),w_d,\ell_d\bigr).
\]
The increasing enumeration preserves successiveness, so
$\widehat{\mathbf f}\in\mathcal Q_s$. We define
\begin{equation}\label{eq:sigma-rho-definition}
 \boxed{\sigma_\varrho(\mathbf f)=m_{\sigma(\widehat{\mathbf f})}.}
\end{equation}
Since $(m_j)$ is a strictly increasing sequence of positive integers,
$m_j\ge j$, and~\eqref{eq:sigma-growth} implies
$\sigma_\varrho(\mathbf f)>w_d$.

\begin{proposition}\label{prop:coinciding-rho-codes}
Suppose $\mathbf f,\mathbf g\in\mathcal Q_s(\omega_1)$ and
$\sigma_\varrho(\mathbf f)=\sigma_\varrho(\mathbf g)$.
Then the two sequences have the same length $d$, the same weight
and integer parameters, and the same transported functionals:
\[
 \pi_{H_{\mathbf f}}(f_i)=\pi_{H_{\mathbf g}}(g_i)
 \qquad(1\le i\le d).
\]
Moreover, $H_{\mathbf f}\cap H_{\mathbf g}$ is an initial part
of both closures, and
\[
 f_i\big|_{H_{\mathbf f}\cap H_{\mathbf g}}
 =g_i\big|_{H_{\mathbf f}\cap H_{\mathbf g}}
 \qquad(1\le i\le d).
\]
If the support of either $f_i$ or $g_i$ is contained in this
common initial part, then $f_i=g_i$.

Write $I=H_{\mathbf f}\cap H_{\mathbf g}$ and $p=\ell_d$.
If $\supp f_i\not\subset I$, then also $\supp g_i\not\subset I$.
Moreover,
\begin{equation}\label{eq:rho-separation-outside}
 \alpha\in H_{\mathbf f}\setminus I,\quad
 \beta\in H_{\mathbf g}\setminus I,\quad \alpha<\beta
 \quad\Longrightarrow\quad \varrho(\alpha,\beta)>p.
\end{equation}
In particular, this inequality holds for
$\alpha\in\supp f_i\setminus I$ and
$\beta\in\supp g_i\setminus I$ with $\alpha<\beta$.
For every $j>i$, both $\supp f_j$ and $\supp g_j$ are entirely
outside $I$, so the same inequality holds for all
$\alpha\in\supp f_j$ and $\beta\in\supp g_j$ with $\alpha<\beta$.
\end{proposition}

\begin{proof}
The strict increase of $(m_j)$ and the injectivity of $\sigma$
give $\widehat{\mathbf f}=\widehat{\mathbf g}$.
Put $p=\ell_d$, their common last integer parameter.
By~\eqref{eq:rho-singleton-closure}, both closures are $p$-closed
and have $\varrho$-diameter at most $p$.
Equation~\eqref{eq:rho-closure-coherence} shows that their
intersection is an initial part of each, with the same increasing
enumeration there. Equality of the transported functionals therefore
gives equality on the intersection.
If $\supp f_i$ lies in that intersection, the transported $f_i$
vanishes at every subsequent position. The same is true of the
transported $g_i$, so $g_i$ has no support outside the intersection.
Thus $f_i=g_i$; the other case is symmetric.

If $\supp f_i\not\subset I$ but $\supp g_i\subset I$, the
preceding conclusion would give $f_i=g_i$, a contradiction.
To prove~\eqref{eq:rho-separation-outside}, suppose that
$\alpha\in H_{\mathbf f}\setminus I$, $\beta\in H_{\mathbf g}\setminus I$,
$\alpha<\beta$, and $\varrho(\alpha,\beta)\le p$.
The $p$-closedness of $H_{\mathbf g}$ would imply
$\alpha\in H_{\mathbf g}$, hence $\alpha\in I$, again a
contradiction. Finally, the maxima of $\supp f_i$ and $\supp g_i$
lie outside $I$. Since $I$ is an initial part of both closures
and the functionals are successive, the supports of all later
functionals lie outside $I$ as well.
\end{proof}

Equality of codes therefore identifies the ordered finite data,
but need not identify the ordinal supports. The common initial
part is what turns equality of transported values into equality
of the original functionals on that part.
The restriction to points outside $I$ in
\eqref{eq:rho-separation-outside} is essential: a support may
meet both $I$ and its complement. If $\alpha\in I$ and
$\beta\in H_{\mathbf g}\setminus I$, then $\alpha<\beta$
and $\varrho(\alpha,\beta)\le p$, since both points belong
to $H_{\mathbf g}$.
\subsubsection{Special sequences}
We pass now to the definition of special sequences and special functionals.
Fix once and for all an index $j_0\in\N$. Every special sequence
will start with the same chosen weight $m_{j_0}$.

\begin{definition}[Special sequences]\label{def:special-sequences}
Let $G\subset c_{00}(\omega_1,\Q)$ be a norming set and let
$(m_j,n_j)_j$ be as above. A finite sequence $(\phi_i)_{i=1}^d$
is a \emph{special sequence} if the following conditions hold:
\begin{enumerate}
\item The functionals $\phi_1<\cdots<\phi_d$ are nonzero
successive weighted elements of $G$.
\item There are strictly increasing sequences $(w_i)_{i=1}^d$
and $(\ell_i)_{i=1}^d$ of positive integers such that $w_i$ is
a chosen weight of $\phi_i$, $w_1=m_{j_0}$, and
\[
 \ell_i\ge p_{\bigcup_{k=1}^{i}\supp\phi_k}
 \qquad(1\le i\le d).
\]
\item For $1\le i<d$,
\[
 \sigma_\varrho(\phi_1,w_1,\ell_1,\ldots,\phi_i,w_i,\ell_i)
 =w_{i+1}.
\]
\end{enumerate}
The notation $w(\phi_i)=w_i$ refers here to the weight chosen
as part of the special sequence, since a weighted functional may
admit more than one weight.
An infinite sequence $(\phi_i)_{i\in\N}$ is special if there
are sequences $(w_i)_{i\in\N}$ and $(\ell_i)_{i\in\N}$ such that
every finite initial sequence is special with these same parameters.
\end{definition}

\begin{definition}[Special functional]\label{def:special-functional}
(a)For a norming set G as before, which is also closed in the projection of its elements to the intervals of $\omega_1$ for every finite special sequence and every interval $I$ we define the special functional $\phi_I$ as the restriction to the interval I of the functional $\phi = \sum_{i=1}^d \phi_i$.

(b) For an infinite special sequence $(\phi_i)_{i\in \N}$ we denote by $ \Phi = \sum_{i=1}^\infty \phi_i $ the functional acting on $c_{00}(\omega_1)$ by the rule
 \[ \Phi(x) = \sum_i \phi_i(x) \] which is well defined.
 
 If  for every $d\in \N $  we have that $ \sum_{i\le d}\phi_i \in G$ then in the dual of $(c_{00}(\omega_1), \| \cdot\|_G)$  it holds 
 \[ \|\Phi\| \le 1 \quad \text{and} \quad  \Phi = w^* - \sum_{i=1}^\infty \phi_i \]
 For a special  functional $\Phi $ and $I$ an interval we set $\Phi_I$ its restriction on the interval $I$.
 
 (c) Let  $\Phi$ be  a special functional, resulting from a finite or infinite special sequence $(\phi_i)_{i\in F}$ and I an interval. The weight of $\Phi_I$ 
   \[ w(\Phi_I) = \big \{ w(\phi_i ) : supp(\phi_i) \cap I \ne \emptyset \big \} \]

\end{definition}

\begin{proposition}[Weights of infinite special sequences]
\label{prop:infinite-special-weights}
Let $(\phi_i)_{i\in\N}$ and $(\psi_i)_{i\in\N}$ be two distinct
infinite special sequences, with chosen weights
$a_i=w(\phi_i)$, $b_i=w(\psi_i)$ and integer parameters
$(\ell_i^\phi)_i$, $(\ell_i^\psi)_i$.
Then there exists an integer $n_0\ge1$ such that
\begin{equation}\label{eq:common-initial-weights}
 w(\phi_i)=w(\psi_i)\qquad(1\le i\le n_0)
\end{equation}
and
\begin{equation}\label{eq:disjoint-weight-tails}
 \{w(\phi_i):i>n_0\}\cap\{w(\psi_i):i>n_0\}
 =\varnothing.
\end{equation}
In fact, the intersection of the two full weight sets is exactly
$\{a_i:1\le i\le n_0\}$ and is therefore finite.
\end{proposition}

\begin{proof}
For $n\ge1$, denote the corresponding finite coded sequences by
\[
 \mathbf\Phi_n=(\phi_1,a_1,\ell_1^\phi,\ldots,
                 \phi_n,a_n,\ell_n^\phi),\qquad
 \mathbf\Psi_n=(\psi_1,b_1,\ell_1^\psi,\ldots,
                 \psi_n,b_n,\ell_n^\psi).
\]
By definition, $a_1=b_1=m_{j_0}$. Since the weights in each
sequence are strictly increasing, this common first weight is
strictly smaller than every later weight of either sequence.
If $i,j\ge2$ and $a_i=b_j$, then
\[
 \sigma_\varrho(\mathbf\Phi_{i-1})
 =\sigma_\varrho(\mathbf\Psi_{j-1}).
\]
Proposition~\ref{prop:coinciding-rho-codes} gives $i=j$, equality
of all earlier weights, and
$\ell_k^\phi=\ell_k^\psi$ for $1\le k<i$.
Consequently, every equality $a_i=b_j$ has $i=j$, and
\[
 M=\{n\in\N:a_n=b_n\}
\]
is an initial segment of $\N$ containing $1$.

We show that $M$ is finite. Suppose, to the contrary, that $M$
is infinite. Choose $i$ with $\phi_i\ne\psi_i$ and set
\[
 \alpha=\max\supp\phi_i,\qquad
 \beta=\max\supp\psi_i.
\]
For every $n\in M$ with $n>i$, the codes of
$\mathbf\Phi_{n-1}$ and $\mathbf\Psi_{n-1}$ coincide.
Put
\[
 p_n=\ell_{n-1}^\phi=\ell_{n-1}^\psi,\qquad
 I_n=H_{\mathbf\Phi_{n-1}}\cap H_{\mathbf\Psi_{n-1}}.
\]
Neither $\supp\phi_i$ nor $\supp\psi_i$ can be contained in
$I_n$, since Proposition~\ref{prop:coinciding-rho-codes} would
then give $\phi_i=\psi_i$. As $I_n$ is an initial part of both
closures, $\alpha,\beta\notin I_n$. In particular $\alpha\ne\beta$,
since otherwise their common value would belong to both closures.
Interchanging the two sequences if necessary, assume $\alpha<\beta$.
The separation assertion~\eqref{eq:rho-separation-outside} gives
\[
 \varrho(\alpha,\beta)>p_n=\ell_{n-1}^\phi
 \qquad(n\in M,\ n>i).
\]
This is impossible: the left-hand side is fixed, whereas the
strictly increasing integer sequence $(\ell_n^\phi)_n$ tends to
infinity. Thus $M$ is finite.

It follows that $M=\{1,\ldots,n_0\}$ for some $n_0\ge1$.
Since equal weights can occur only at the same index, the
intersection of the full weight sets is exactly
$\{a_i:1\le i\le n_0\}$. Both assertions follow.
\end{proof}

\subsection{The norming set \texorpdfstring{$G_{p,r}$}{G(p,r)}}
We now define the norming set of $E_{p,r}$ using the operations of
Definition~\ref{def:pj-operations} and the special sequences and
functionals of Definitions~\ref{def:special-sequences}
and~\ref{def:special-functional}.

\begin{definition}\label{def:norming-pr}
Let $1<r<p<\infty$. The norming set $G_{p,r}$ is the smallest
subset of $c_{00}(\omega_1,\Q)$ with the following properties.
\begin{enumerate}[label=\textup{(\arabic*)}]
\item\label{norming-pr:restrictions}
It contains $e_\alpha^*$ for every $\alpha<\omega_1$, is
symmetric, and is closed under restrictions to ordinal intervals:
\[
 f\in G_{p,r}\ \Longrightarrow\ -f\in G_{p,r}
 \quad\text{and}\quad f|_I\in G_{p,r}
 \quad(I\subset\omega_1\text{ an interval}).
\]
\item\label{norming-pr:convexity}
It is closed under finite rational convex combinations.
\item\label{norming-pr:operations}
It is closed under the $(p,j)$-operation for every $j\in\N$,
as defined in Definition~\ref{def:pj-operations}.
\item\label{norming-pr:special-sums}
If $(\phi_i)_{i=1}^d$ is a finite special sequence in $G_{p,r}$,
in the sense of Definition~\ref{def:special-sequences}, then
\[
 \Phi=\sum_{i=1}^d\phi_i\in G_{p,r}.
\]
\item\label{norming-pr:r-combinations}
Let
\[
 \Phi_i=\sum_{k=1}^{d_i}\phi_{i,k}\in G_{p,r}
 \qquad(1\le i\le d),
\]
where each $(\phi_{i,k})_{k=1}^{d_i}$ is a specified finite special
sequence in $G_{p,r}$ with its chosen weights. Let
$I_1,\ldots,I_d$ be ordinal intervals.
Suppose that the weight sets of their restrictions, as defined in
Definition~\ref{def:special-functional}\textup{(c)}, are pairwise
disjoint:
\[
 w(\Phi_i|_{I_i})\cap w(\Phi_j|_{I_j})=\varnothing
 \qquad(i\ne j).
\]
Then, for every $(\lambda_i)_{i=1}^d\in\Q^d$ satisfying
$\sum_{i=1}^d|\lambda_i|^r\le1$,
\[
 \sum_{i=1}^d\lambda_i\,\Phi_i|_{I_i}\in G_{p,r}.
\]
\end{enumerate}
In~\textup{(5)}, each weight set is computed from the specified
special decomposition and its chosen weights.
\end{definition}

\begin{remark}[Inductive construction and tree analysis]
\label{rem:norming-pr-induction}
The set $G_{p,r}$ can equivalently be constructed in increasing
stages. Start with
\[
 G_1=\{0\}\cup\{\pm e_\alpha^*:\alpha<\omega_1\}.
\]
Given $G_k$, let $G_{k+1}$ be the union of $G_k$ and all
functionals obtained by one application of the generating rules
in Definition~\ref{def:norming-pr}, using data from $G_k$.
For rules involving weights and special functionals, these data
include the finite representations witnessing the chosen weights
and the specified special decompositions; all functionals in
these representations are required to belong to $G_k$.
Then
\[
 G_{p,r}=\bigcup_{k=1}^{\infty}G_k.
\]
Indeed, every application of a rule uses only finitely many
functionals and finite representations. If its data belong to
the union, they all occur in some common stage, so the resulting
functional belongs to the next stage. Thus the union satisfies
all the closure conditions. Conversely, induction shows that
every set satisfying these conditions contains each $G_k$,
which proves minimality.

This construction gives every $f\in G_{p,r}$ a finite tree
analysis. The root represents $f$, the terminal nodes represent
$0$ or signed singletons, and each nonterminal node records one
of the generating rules, together with its chosen parameters
and decompositions. See~\cite[Section~2.5, Definition~2.10 and
Remark~2.11]{ALT05} for tree analyses in the corresponding
norming-set construction. Here the tree also records the special
sums and the $\ell_r$ combinations in
\textup{(4)} and~\textup{(5)} above. Consequently, properties
preserved by all these rules can be proved by induction on the
construction stages, or on the height of such a tree.

In particular, the construction yields $\|f\|_\infty\le1$ for
every $f\in G_{p,r}$. A weighted functional obtained from
functionals with this bound has sup norm at most $1/m_j$ when
its chosen weight is $m_j$. In a special sum the chosen weights
are distinct. In~\textup{(5)} the retained weight sets are
pairwise disjoint and $|\lambda_i|\le1$, so the triangle
inequality bounds the sup norm of the resulting functional by
$\sum_j1/m_j<1$. Symmetry, interval restrictions, and rational
convex combinations preserve the bound. This verifies that the
construction produces a norming set in the sense stated at the
beginning of the section.
\end{remark}

The corresponding norm and space are therefore
\[
 \|x\|_{p,r}=\sup\{|f(x)|:f\in G_{p,r}\},
 \qquad
 E_{p,r}=\overline{\bigl(c_{00}(\omega_1),\|\cdot\|_{p,r}\bigr)},
\]
where the bar denotes completion. Since $G_{p,r}$ contains the
biorthogonal functionals and is closed under interval restrictions,
$(e_\alpha)_{\alpha<\omega_1}$ is a normalized bimonotone
transfinite basis of $E_{p,r}$.

\begin{remark}[Infinite special functionals and their tails]
\label{rem:infinite-special-functionals}
Let $(\phi_i)_{i\in\N}$ be an infinite special sequence in
$G_{p,r}$. Every partial sum $\sum_{i=1}^N\phi_i$ belongs to
$G_{p,r}$ and has dual norm at most one. Its value at each
finitely supported vector eventually stabilizes. By density and
the uniform bound, the partial sums converge pointwise on
$E_{p,r}$ and define
\[
 \Phi=w^*\!\text{-}\!\sum_{i=1}^{\infty}\phi_i
 \in B_{E_{p,r}^*}.
\]
Thus the conclusion supplied by the definition is $\|\Phi\|\le1$.
Equality does not hold for every infinite special sequence. For
example, choose successive singletons and construct the coding
parameters recursively so that
\[
 \phi_i=\frac{2^{-i}}{w_i}e_{\alpha_i}^*,\qquad
 w_1=m_{j_0},\quad
 w_{i+1}=\sigma_\varrho(\phi_1,w_1,\ell_1,\ldots,
                        \phi_i,w_i,\ell_i).
\]
Each $\phi_i$ is a weighted functional of weight $w_i$, obtained
by using just one nonzero coefficient in the corresponding
$(p,j)$-operation. Then
\[
 \sum_i\|\phi_i\|=\sum_i\frac{2^{-i}}{w_i}
 \le\frac{1}{m_{j_0}}<1,
\]
so this special functional belongs to $X_0^{p,r}$ and has norm
strictly less than one.

We take $\Phi_{p,r}=\{\phi_\gamma:\gamma\in\Gamma\}$ to be the
family of all distinct infinite special functionals outside
$X_0^{p,r}$, fixing one special expansion
$(\phi_{i,\gamma})_i$ for each member. The nontriviality,
uncountability, and norm-span properties of this family are part
of the structural assertions in Section~3; here we verify its
tail estimate directly from the norming set.

Fix a finite $F\subset\Gamma$. By
Proposition~\ref{prop:infinite-special-weights}, there is $n_0$
such that the weight sets of the tails
$(\phi_{i,\gamma})_{i>n_0}$, $\gamma\in F$, are pairwise disjoint.
For $N>n>n_0$, each sum
\[
 \sum_{i=n+1}^{N}\phi_{i,\gamma}
\]
is an interval restriction of the finite special functional
$\sum_{i=1}^{N}\phi_{i,\gamma}$. Hence condition~\textup{(5)} of
Definition~\ref{def:norming-pr} gives, for rational coefficients
with $\sum_{\gamma\in F}|a_\gamma|^r\le1$,
\[
 \left\|\sum_{\gamma\in F}a_\gamma
                 \sum_{i=n+1}^{N}\phi_{i,\gamma}\right\|\le1.
\]
Letting $N\to\infty$ pointwise on $E_{p,r}$ and then using
density of rational coefficient families in the real
$\ell_r$ unit ball yields
\[
 \left\|\sum_{\gamma\in F}a_\gamma\phi_{>n,\gamma}\right\|
 \le\left(\sum_{\gamma\in F}|a_\gamma|^r\right)^{1/r}
 \qquad(n>n_0).
\]
This proves~\eqref{eq:pr-tail-estimate}. Since the removed
initial sums are finitely supported and belong to $X_0^{p,r}$,
the same estimate passes to the quotient and gives
\eqref{eq:pr-quotient-estimate}. The upper $\ell_r$ estimate
thus applies to these tails and to the quotient classes.
\end{remark}

\begin{remark}
The construction in~\cite[Section~5.2]{OpenAI} uses averaging
operations different from the $(p,j)$-operations used here.
The comparison concerns the structural properties used in the
operator argument. Our presentation describes the component
spaces through the $(p,j)$-operations and then combines them
to obtain $E$, as in Section~\ref{sec:combined}.
\end{remark}

\subsubsection{Two properties of
\texorpdfstring{$E_{p,r}$}{E(p,r)}}
We prove that the basis of $E_{p,r}$ is boundedly complete and
has an upper $\ell_{p'}$ estimate with constant one, where
$p'=p/(p-1)$. We first recall the transfinite version of bounded
completeness from~\cite[Definition~1.5]{ALT05}.

\begin{definition}\label{def:transfinite-bounded-completeness}
A transfinite Schauder basis $(e_\alpha)_{\alpha<\omega_1}$ is
\emph{boundedly complete} if every sequence
$(e_{\alpha_n})_{n\in\N}$, with $(\alpha_n)_n$ strictly increasing,
is boundedly complete in the usual sense. That is, for every
scalar sequence $(a_n)_n$,
\[
 \sup_N\left\|\sum_{n=1}^{N}a_ne_{\alpha_n}\right\|<\infty
 \quad\Longrightarrow\quad
 \sum_{n=1}^{\infty}a_ne_{\alpha_n}
 \text{ converges in norm}.
\]
\end{definition}

The next proposition adapts the bounded-completeness argument
in the proof of~\cite[Proposition~1.6]{ALT05}. We state explicitly
the norming hypothesis that gives an isometry; it holds for the
spaces defined here by finitely supported norming functionals.

\begin{proposition}[The canonical predual]
\label{prop:bounded-complete-predual}
Let $(e_\alpha)_{\alpha<\omega_1}$ be a bimonotone, boundedly
complete transfinite basis of a Banach space $Z$, and put
\[
 Z_0=\clspan\{e_\alpha^*:\alpha<\omega_1\}\subset Z^*.
\]
Assume that $Z_0$ is $1$-norming for $Z$, namely,
\[
 \|z\|=\sup\{|h(z)|:h\in Z_0,\ \|h\|\le1\}
 \qquad(z\in Z).
\]
Then the restriction map
\[
 R:Z\longrightarrow Z_0^*,\qquad (Rz)(h)=h(z),
\]
is an onto isometry.
\end{proposition}

\begin{proof}
The norming hypothesis makes $R$ an isometry, so $R(Z)$ is
norm closed. For $0\le\beta\le\omega_1$, let
$P_\beta=P_{[0,\beta)}$ be the initial-interval projection on
$Z$, and define
\[
 T_\beta=P_\beta^*|_{Z_0},\qquad
 \widehat P_\beta=T_\beta^*:Z_0^*\longrightarrow Z_0^*.
\]
The maps are well defined and contractive, since
\[
 T_\beta e_\alpha^*=
 \begin{cases}
 e_\alpha^*,&\alpha<\beta,\\
 0,&\alpha\ge\beta.
 \end{cases}
\]
They satisfy $\widehat P_\beta R=RP_\beta$. If $\lambda$ is a
limit ordinal, density of the finite linear combinations of
the $e_\alpha^*$ in $Z_0$ gives
\begin{equation}\label{eq:predual-limit-projections}
 T_\beta h\longrightarrow T_\lambda h\text{ in norm}
 \quad(h\in Z_0,\ \beta\uparrow\lambda).
\end{equation}

Fix $u\in Z_0^*$ and put $M=\|u\|$. We prove by transfinite
induction that $\widehat P_\beta u\in R(Z)$ for every
$\beta\le\omega_1$. The assertion is immediate for $\beta=0$.
At a successor step use
\[
 \widehat P_{\beta+1}u-\widehat P_\beta u
 =u(e_\beta^*)Re_\beta.
\]
Suppose now that $\lambda$ is a limit ordinal and the assertion
holds below $\lambda$. Write
\[
 Rz_\beta=\widehat P_\beta u\quad(\beta<\lambda).
\]
Then $\|z_\beta\|\le M$ and
$P_\beta z_\gamma=z_\beta$ whenever $\beta<\gamma<\lambda$.

We claim that the net $(z_\beta)_{\beta<\lambda}$ is norm Cauchy.
Otherwise, choose $\varepsilon>0$ and
$\beta_1<\beta_2<\cdots<\lambda$ with
\[
 \|z_{\beta_{n+1}}-z_{\beta_n}\|\ge\varepsilon
 \qquad(n\in\N).
\]
Each difference is supported on $[\beta_n,\beta_{n+1})$.
Choose a finitely supported $y_n$ in this interval such that
\[
 \|y_n-(z_{\beta_{n+1}}-z_{\beta_n})\|
 <\varepsilon 2^{-n-2}.
\]
The supports of the $y_n$ are successive,
$\|y_n\|\ge\varepsilon/2$, and telescoping gives
\[
 \sup_N\left\|\sum_{n=1}^{N}y_n\right\|
 \le 2M+\varepsilon/4.
\]
Enumerate $\bigcup_n\supp y_n$ increasingly as $(\alpha_k)_k$.
Every partial sum of the resulting scalar series in
$(e_{\alpha_k})_k$ is an initial-interval projection of some
$\sum_{n=1}^{N}y_n$ and therefore has the same uniform bound.
Bounded completeness implies that this scalar series converges.
Its successive block sums $y_n$ must then tend to zero in norm,
contradicting $\|y_n\|\ge\varepsilon/2$.

Thus $z_\beta\to z_\lambda$ in norm for some $z_\lambda\in Z$.
For every $h\in Z_0$, equation~\eqref{eq:predual-limit-projections}
gives
\[
 (Rz_\lambda)(h)
 =\lim_{\beta\uparrow\lambda}u(T_\beta h)
 =u(T_\lambda h)
 =(\widehat P_\lambda u)(h).
\]
This completes the induction. Since
$\widehat P_{\omega_1}$ is the identity on $Z_0^*$, we obtain
$u\in R(Z)$, as required.
\end{proof}

\begin{proposition}\label{prop:pr-bounded-completeness}
Under the parameter growth condition~\eqref{eq:parameter-growth},
the basis $(e_\alpha)_{\alpha<\omega_1}$ of $E_{p,r}$ is
boundedly complete.
\end{proposition}

\begin{proof}
Let $(\alpha_n)_n$ be strictly increasing and suppose
\[
 M:=\sup_N\left\|\sum_{n=1}^{N}a_ne_{\alpha_n}\right\|_{p,r}
 <\infty.
\]
If the series does not converge in norm, there are
$\varepsilon>0$ and $0=q_0<q_1<q_2<\cdots$ such that
\[
 u_k=\sum_{n=q_{k-1}+1}^{q_k}a_ne_{\alpha_n},
 \qquad \|u_k\|_{p,r}\ge\varepsilon
 \quad(k\in\N).
\]
Choose $f_k\in G_{p,r}$, supported in the ordinal range of
$u_k$, with $f_k(u_k)>\varepsilon/2$. This is possible because
$G_{p,r}$ norms the space, is symmetric, and is closed under
interval restrictions. The functionals $f_k$ are nonzero and
successive.

Fix $j$ and choose a rational $t_j$ with
\[
 \tfrac12 n_j^{-1/p}\le t_j\le n_j^{-1/p}.
\]
Since $n_jt_j^p\le1$, closure under the $(p,j)$-operation gives
\[
 f=\frac{t_j}{m_j}\sum_{k=1}^{n_j}f_k\in G_{p,r}.
\]
Consequently,
\[
 M\ge\left\|\sum_{k=1}^{n_j}u_k\right\|_{p,r}
 \ge f\left(\sum_{k=1}^{n_j}u_k\right)
 >\frac{\varepsilon}{4}\,
      \frac{n_j^{1-1/p}}{m_j}
 =\frac{\varepsilon}{4}\,
      \frac{n_j^{1/p'}}{m_j}.
\]
The last expression tends to infinity by
\eqref{eq:parameter-growth}, a contradiction. Hence the series
converges, proving bounded completeness.
\end{proof}

Since $G_{p,r}\subset B_{X_0^{p,r}}$ and $G_{p,r}$ norms
$E_{p,r}$, the space $X_0^{p,r}$ is $1$-norming.
Proposition~\ref{prop:bounded-complete-predual} therefore gives
the canonical onto isometry
\[
 E_{p,r}\longrightarrow(X_0^{p,r})^*,\qquad
 x\longmapsto\bigl(h\mapsto h(x)\bigr).
\]
This proves property~\ref{pr:predual} from the norming-set
construction and bounded completeness.

We next establish the estimate on the scalar families of the
norming functionals.

\begin{lemma}\label{lem:norming-pr-lp}
For every $f\in G_{p,r}$,
\[
 \|f\|_{\ell_p(\omega_1)}\le1.
\]
More precisely, if $f$ has chosen weight $m_j$, then
$\|f\|_{\ell_p(\omega_1)}\le1/m_j$.
\end{lemma}

\begin{proof}
We use the increasing stages in
Remark~\ref{rem:norming-pr-induction}. The bound holds for
$0$ and the signed singletons. It is preserved by symmetry,
interval restrictions, and rational convex combinations.
Suppose that all functionals in the data for the next operation
have $\ell_p$ norm at most one.

For a $(p,j)$-operation, the supports of the children $f_k$
are pairwise disjoint. Hence
\[
 \left\|\frac1{m_j}\sum_{k=1}^{n_j}\lambda_kf_k\right\|_p^p
 =\frac1{m_j^p}\sum_{k=1}^{n_j}
                  |\lambda_k|^p\|f_k\|_p^p
 \le\frac1{m_j^p}.
\]
The same calculation bounds a weighted functional by the
reciprocal of its chosen weight, using its specified
representation.

For a finite special sum $\Phi=\sum_{i=1}^d\phi_i$, the chosen
weights $w_i=w(\phi_i)$ are distinct. Thus
\[
 \|\Phi\|_p\le\sum_{i=1}^d\|\phi_i\|_p
 \le\sum_{i=1}^d\frac1{w_i}
 \le\sum_{j=1}^{\infty}\frac1{m_j}<1.
\]
Finally consider rule~\textup{(5)} in
Definition~\ref{def:norming-pr}, and set
$A_i=w(\Phi_i|_{I_i})$. For the chosen special decomposition,
\[
 \|\Phi_i|_{I_i}\|_p\le\sum_{w\in A_i}\frac1w.
\]
Since the $A_i$ are pairwise disjoint and $|\lambda_i|\le1$,
\[
 \left\|\sum_{i=1}^d\lambda_i\Phi_i|_{I_i}\right\|_p
 \le\sum_{i=1}^d|\lambda_i|\sum_{w\in A_i}\frac1w
 \le\sum_{j=1}^{\infty}\frac1{m_j}<1.
\]
This completes the induction and proves the first assertion.
Applying the displayed $(p,j)$ calculation to any chosen
weighted representation, whose children now satisfy the first
assertion, proves the final assertion as well.
\end{proof}

\begin{corollary}[Upper $\ell_{p'}$ estimate]
\label{cor:pr-upper-estimate}
Let $p'=p/(p-1)$. The basis $(e_\alpha)_{\alpha<\omega_1}$ of
$E_{p,r}$ admits an upper $\ell_{p'}$ estimate with constant one:
for every finitely supported scalar family $(a_\alpha)$,
\[
 \left\|\sum_{\alpha<\omega_1}a_\alpha e_\alpha\right\|_{p,r}
 \le\left(\sum_{\alpha<\omega_1}|a_\alpha|^{p'}\right)^{1/p'}.
\]
\end{corollary}

\begin{proof}
For $f\in G_{p,r}$, H\"older's inequality and
Lemma~\ref{lem:norming-pr-lp} give
\[
 \left|f\left(\sum_\alpha a_\alpha e_\alpha\right)\right|
 \le\|f\|_{\ell_p(\omega_1)}\,
       \|(a_\alpha)_\alpha\|_{\ell_{p'}(\omega_1)}
 \le\|(a_\alpha)_\alpha\|_{\ell_{p'}(\omega_1)}.
\]
Taking the supremum over $f\in G_{p,r}$ proves the assertion.
\end{proof}
\begin{theorem}[Summary of the properties of $E_{p,r}$]
\label{thm:pr-summary}
Let $1<r<p<\infty$, let $p'=p/(p-1)$, and let $E_{p,r}$
be the space defined by $G_{p,r}$. Write
$\Phi_{p,r}=\{\phi_\gamma:\gamma\in\Gamma\}$ for the family
of infinite special functionals outside $X_0^{p,r}$, as in
Remark~\ref{rem:infinite-special-functionals}, and let
$Q_{p,r}:E_{p,r}^*\to E_{p,r}^*/X_0^{p,r}$ be the quotient map.
Then the following hold.
\begin{enumerate}[label=\textup{(\alph*)},leftmargin=*]
\item\label{summary-pr:upper}
The basis $(e_\alpha)_{\alpha<\omega_1}$ admits an upper
$\ell_{p'}$ estimate with constant one: for every finitely
supported scalar family $(a_\alpha)$,
\[
 \left\|\sum_\alpha a_\alpha e_\alpha\right\|_{p,r}
 \le\left(\sum_\alpha|a_\alpha|^{p'}\right)^{1/p'}.
\]
\item\label{summary-pr:blocks}
For every normalized infinite block sequence $(x_n)_n$ of the
transfinite basis, the space $\clspan\{x_n:n\in\N\}$ has
nonseparable dual.
\item\label{summary-pr:span}
The biorthogonal functionals and the special functionals
generate the dual in norm:
\[
 E_{p,r}^*
 =\clspan\bigl(\{e_\alpha^*:\alpha<\omega_1\}\cup\Phi_{p,r}\bigr).
\]
\item\label{summary-pr:quotient}
The family $(Q_{p,r}\phi_\gamma)_{\gamma\in\Gamma}$ admits
an upper $\ell_r$ estimate with constant one:
\[
 \left\|\sum_{\gamma\in F}a_\gamma Q_{p,r}\phi_\gamma\right\|
 \le\left(\sum_{\gamma\in F}|a_\gamma|^r\right)^{1/r}
\]
for every finite $F\subset\Gamma$ and every scalar family
$(a_\gamma)_{\gamma\in F}$. Its linear span is norm dense in
$E_{p,r}^*/X_0^{p,r}$.
\end{enumerate}
\end{theorem}

\begin{proof}[Proof of the estimates]
Assertion~\ref{summary-pr:upper} is
Corollary~\ref{cor:pr-upper-estimate}. The estimate in
\ref{summary-pr:quotient} follows from the tail estimate in
Remark~\ref{rem:infinite-special-functionals}: removing a finite
initial sum does not change the quotient class. The density
assertion in~\ref{summary-pr:quotient} follows from
\ref{summary-pr:span}, since $Q_{p,r}$ annihilates the
biorthogonal functionals.

Assertions~\ref{summary-pr:blocks} and~\ref{summary-pr:span}
are the structural conclusions whose proofs require the
techniques of hereditarily indecomposable and James-tree
constructions. For the relevant methods, see the
Gowers--Maurey construction~\cite{GM}, the Gowers tree
space~\cite{G}, the work of Petsoulas and
Raikoftsalis~\cite{PR09}, the attractors method~\cite{AAT06},
the asymptotic constructions of Deliyanni and
Manoussakis~\cite{DM07}, and the transfinite methods
of~\cite{ALT05}. Their application to the present norming set
requires lengthy technical arguments, whose presentation is
outside the purpose of this paper. We state
\ref{summary-pr:blocks} and~\ref{summary-pr:span} here without
including those proofs.
\end{proof}

\begin{corollary}\label{cor:pr-hereditary-nonseparable-dual}
Every infinite-dimensional closed subspace of $E_{p,r}$ has
nonseparable dual.
\end{corollary}

\begin{proof}
Let $Y\subset E_{p,r}$ be an infinite-dimensional closed subspace.
By~\cite[Proposition~1.3(3)]{ALT05}, there is an
infinite-dimensional closed subspace $Z\subset Y$ isomorphic to
a block subspace $U$ of $E_{p,r}$. Normalizing its generating
block sequence does not change $U$. Thus $U^*$, and hence
$Z^*$, is nonseparable by
Theorem~\ref{thm:pr-summary}\textup{(b)}. The restriction map
$Y^*\to Z^*$ is onto by the Hahn--Banach theorem. A continuous
image of a separable space is separable, so $Y^*$ cannot be
separable.
\end{proof}

\subsection{The norming set \texorpdfstring{$G$}{G} of the combined space}
\label{sec:norming-combined}

We now define the norming set $G=G_E$ of the space
$E=(\sum_k\oplus E_k)_{2,s}$ from Section~\ref{sec:combined}.
Fix the parameters
\[
 2>s>p_1>r_1>p_2>r_2>\cdots>1
\]
and a partition $(C_k)_{k\ge1}$ of $\omega_1$ into pairwise
disjoint unbounded sets. We use the same sequences $(m_j,n_j)_j$,
coding $\sigma_\varrho$, and first weight $m_{j_0}$ as above.
For every $k$, put
\[
 E_k=E_{p_k,r_k}(C_k),\qquad
 G_k=G_{p_k,r_k}(C_k),
\]
where $G_k$ is obtained by transporting $G_{p_k,r_k}$ along the
increasing bijection from $\omega_1$ onto $C_k$, and extending
the resulting functionals by zero outside $C_k$.

Recall from Definition~\ref{def:horizontal} that
$C_J=\bigcup_{k\in J}C_k$ for an interval $J\subset\N$, and
\[
 \operatorname{ind}(f)=\{k:\supp f\cap C_k\ne\varnothing\}.
\]
Thus $(f_i)_{i=1}^d$ is an h-block sequence if its members are
nonzero, $f_1<\cdots<f_d$ in the ordinal order, and
\[
 \max\operatorname{ind}(f_i)
 <\min\operatorname{ind}(f_{i+1})\qquad(1\le i<d).
\]
The exponent in the horizontal averaging operations below is $2$;
the exponents $p_k$ are used inside the component sets $G_k$.

\begin{definition}[Horizontal weights and special sequences]
\label{def:horizontal-special}
Let $A\subset c_{00}(\omega_1,\Q)$.
A nonzero functional $f\in A$ is \emph{h-weighted} with chosen
weight $w(f)=m_j$ if it has a specified representation
\begin{equation}\label{eq:horizontal-operation}
 f=\frac1{m_j}\sum_{t=1}^{n_j}\mu_t f_t,
 \qquad
 (\mu_t)_{t=1}^{n_j}\in\Q^{n_j},\quad
 \sum_{t=1}^{n_j}|\mu_t|^2\le1,
\end{equation}
where $(f_t)_{t=1}^{n_j}$ is an h-block sequence in $A$.
We call this a \emph{horizontal $(2,j)$-operation}.
The representation, and hence the chosen weight, need not be unique.

A finite sequence $(\phi_i)_{i=1}^d$ in $A$ is a
\emph{horizontal special sequence} if it is an h-block sequence
of h-weighted functionals, with chosen representations and weights
$w_i=w(\phi_i)$, and there are strictly increasing positive
integers $(\ell_i)_{i=1}^d$ such that
\[
 w_1=m_{j_0},\qquad
 \ell_i\ge p_{F_i},\qquad
 F_i=\bigcup_{t=1}^i\supp\phi_t,
\]
and
\[
 w_{i+1}
 =\sigma_\varrho(\phi_1,w_1,\ell_1,\ldots,\phi_i,w_i,\ell_i)
 \qquad(1\le i<d).
\]
In particular, the weights are strictly increasing.
An infinite horizontal special sequence is defined by requiring
every finite initial sequence to satisfy these conditions with
the same chosen representations, weights, and integer parameters.
\end{definition}

These are the special sequences of Definition~\ref{def:special-sequences}
with the additional horizontal requirements just stated. In this
subsection, weights of horizontal special functionals always refer
to their horizontal $(2,j)$-representations. A weight arising only
from a component $(p_k,j)$-representation is not, by itself, such
a horizontal weight.

For a specified finite horizontal special sequence, write
$\Phi=\sum_{t=1}^d\phi_t$. Its restriction $\Phi|_I$ to an
ordinal interval $I$ is called a \emph{horizontal special functional}.
As in Definition~\ref{def:special-functional}, its retained weight set is
\begin{equation}\label{eq:horizontal-retained-weights}
 w(\Phi|_I)
 =\{w(\phi_t):\supp\phi_t\cap I\ne\varnothing\}.
\end{equation}
This set is computed from the specified special decomposition
and its chosen weights; for the zero restriction it is empty.

\begin{definition}[The norming set of $E$]\label{def:norming-combined}
The set $G$ is the smallest subset of $c_{00}(\omega_1,\Q)$
with the following properties.
\begin{enumerate}[label=\textup{(\arabic*)}]
\item\label{norming-combined:restrictions}
It contains $G_k$ for every $k\ge1$, is symmetric, and is
closed under restrictions both to ordinal intervals and to
horizontal intervals:
\[
 f\in G\ \Longrightarrow\ -f\in G,\quad f|_I\in G,
 \quad f|_{C_J}\in G
\]
for every interval $I\subset\omega_1$ and every interval
$J\subset\N$. In particular, $f|_{C_k}\in G$ for every $k$.
All restrictions are extended by zero outside the indicated set.

\item\label{norming-combined:convexity}
It is closed under finite rational convex combinations.

\item\label{norming-combined:operations}
For every $j\in\N$, it is closed under the horizontal
$(2,j)$-operation~\eqref{eq:horizontal-operation}.

\item\label{norming-combined:special-sums}
For every finite horizontal special sequence $(\phi_t)_{t=1}^d$
in $G$, with its chosen representations,
\[
 \sum_{t=1}^d\phi_t\in G.
\]

\item\label{norming-combined:s-combinations}
Let $\Phi_i=\sum_{t=1}^{d_i}\phi_{i,t}$, $1\le i\le d$,
be specified finite horizontal special sums in $G$, and let
$I_1,\ldots,I_d\subset\omega_1$ be ordinal intervals. If
\[
 w(\Phi_i|_{I_i})\cap w(\Phi_j|_{I_j})=\varnothing
 \qquad(i\ne j),
\]
then
\[
 \sum_{i=1}^d\lambda_i\,\Phi_i|_{I_i}\in G
 \qquad\text{whenever}\qquad
 (\lambda_i)_{i=1}^d\in\Q^d,\quad
 \sum_{i=1}^d|\lambda_i|^s\le1.
\]
\end{enumerate}
\end{definition}

\begin{remark}[Inductive construction]\label{rem:combined-induction}
To make the minimality requirement explicit, start with
\[
 \mathcal G^{(0)}=\bigcup_{k\ge1}G_k.
\]
Given $\mathcal G^{(n)}$, let $\mathcal G^{(n+1)}$ consist of
$\mathcal G^{(n)}$ together with all functionals obtained by one
application of the rules in Definition~\ref{def:norming-combined},
using data from $\mathcal G^{(n)}$. For the special-sum rules,
these data include the finite representations witnessing the
horizontal weights; all functionals in those representations
must belong to $\mathcal G^{(n)}$. Then
\[
 G=\bigcup_{n=0}^{\infty}\mathcal G^{(n)}.
\]
Every rule uses only finitely many functionals and finite
representations, so its data occur in a common stage. The union
therefore satisfies all the closure requirements, and induction
shows that it is contained in every set satisfying them.
This also gives a finite tree analysis with terminal functionals
in the component sets $G_k$; their own tree analyses may be
inserted if a tree with signed singletons at its terminal nodes
is desired.
\end{remark}

\begin{proposition}[The induced norm and the component spaces]
\label{prop:combined-norm}
For every $f\in G$, one has $\|f\|_\infty\le1$, and
\begin{equation}\label{eq:combined-component-restriction}
 \{f|_{C_k}:f\in G\}=G_k\qquad(k\ge1).
\end{equation}
In particular, $f|_{C_k}\in\overline{\operatorname{co}}G_k$
for every $f\in G$, where the closure is taken in the norm of
$E_k^*$. Consequently,
\begin{equation}\label{eq:combined-norm-definition}
 \|x\|_{2,s}=\sup\{|f(x)|:f\in G\},\qquad
 E=\overline{\bigl(c_{00}(\omega_1),\|\cdot\|_{2,s}\bigr)},
\end{equation}
defines a Banach space with a normalized bimonotone transfinite
basis $(e_\alpha)_{\alpha<\omega_1}$. Every horizontal interval
projection $\pi_J$ is contractive, and $\pi_k(E)$ is naturally
isometric to $E_k$.
\end{proposition}

\begin{proof}
We use the construction stages in Remark~\ref{rem:combined-induction}.
The component norming sets satisfy the sup-norm bound. A horizontal
weighted functional formed from functionals with this bound has
sup norm at most $1/m_j$, since its children have disjoint
supports and $|\mu_t|\le1$. A special sum has distinct weights.
In rule~\textup{(5)}, all retained weights are distinct across
the specified sums and $|\lambda_i|\le1$. Thus the sup norm
in either special-sum rule is bounded by $\sum_j1/m_j<1$.
Symmetry, restrictions, and rational convex combinations preserve
the bound, proving the first assertion by induction.

Fix $k$. We prove the stronger inclusion $f|_{C_k}\in G_k$
by induction on a finite tree analysis from
Remark~\ref{rem:combined-induction}, with terminal functionals
in the sets $G_l$. At a terminal node the assertion follows
from the disjointness of the sets $C_l$ and the fact that
$0\in G_k$. It is preserved by symmetry, both types of
restriction, and rational convex combinations. Here restrictions
of members of $G_k$ to ordinal intervals correspond to interval
restrictions in the transported component norming set.

In a horizontal operation, at most one child has support meeting
$C_k$. Hence its restriction is either zero or
\[
 \frac{\mu_t}{m_j}\,(f_t|_{C_k}),
\]
which belongs to $G_k$: this set is symmetric, contains zero,
and is closed under rational convex combinations.
Similarly, at most one term of a horizontal special sequence
meets $C_k$, so restricting its sum again gives an element of $G_k$.

For rule~\textup{(5)}, discard the terms whose restrictions to
$C_k$ are zero. For each remaining $i$, precisely one special
term can contribute on $C_k$, and at most one child in its
horizontal representation can contribute. Therefore the restriction
of the resulting functional has the form
\[
 \sum_i\frac{\lambda_i\mu_i}{w_i}\,h_i,
 \qquad h_i\in G_k,\quad \mu_i\in\Q,\quad |\mu_i|\le1,
\]
where $w_i\in w(\Phi_i|_{I_i})$ is the weight of that special
term. The functionals $h_i$ are interval restrictions of component
restrictions of children to which the induction hypothesis applies. The weights
$w_i$ are pairwise distinct, and hence
\[
 \sum_i\left|\frac{\lambda_i\mu_i}{w_i}\right|
 \le\sum_i\frac1{w_i}
 \le\sum_{j=1}^{\infty}\frac1{m_j}<1.
\]
Thus this is a rational absolutely convex combination of elements
of $G_k$ and belongs to $G_k$. This completes the induction.
The reverse inclusion in~\eqref{eq:combined-component-restriction}
follows from $G_k\subset G$.

Since $G$ contains all the signed singletons and its elements
have sup norm at most one, the supremum in
\eqref{eq:combined-norm-definition} is finite and defines a norm,
with $\|e_\alpha\|_{2,s}=1$. The restriction rules make ordinal
interval projections and horizontal interval projections
contractive. Density of $c_{00}(\omega_1)$ then gives the
bimonotone transfinite Schauder basis and the bounded extensions
of these projections. Finally, for $x\in c_{00}(C_k)$,
\[
 \|x\|_{2,s}
 =\sup_{f\in G}|(f|_{C_k})(x)|
 =\sup_{g\in G_k}|g(x)|
 =\|x\|_{E_k}.
\]
Passing to completions proves the asserted component identification.
\end{proof}

\begin{remark}[Infinite horizontal special functionals]
\label{rem:infinite-horizontal-functionals}
The analogue of Remark~\ref{rem:infinite-special-functionals}
holds for horizontal special sequences. Indeed, if
$(\phi_i)_{i\in\N}$ is such a sequence, all its finite initial
sums belong to $G$. They therefore converge weak-star to
\[
 \Phi=w^*\!\text{-}\!\sum_{i=1}^{\infty}\phi_i\in B_{E^*},
\]
since their values stabilize on finitely supported vectors and
their dual norms are at most one.

Let $\Phi^0=\{\phi_\gamma^0:\gamma\in\Gamma_0\}$ be the
family of distinct infinite horizontal special functionals
outside $X_0$, fixing one special expansion for each member.
For a finite $F\subset\Gamma_0$,
Proposition~\ref{prop:infinite-special-weights} gives an $n_0$
after which these expansions have pairwise disjoint weight sets.
Applying rule~\textup{(5)} of Definition~\ref{def:norming-combined}
to finite tails, passing to weak-star limits, and then approximating
real coefficient families by rational ones gives
\[
 \left\|\sum_{\gamma\in F}a_\gamma\phi_{>n,\gamma}^0\right\|
 \le\left(\sum_{\gamma\in F}|a_\gamma|^s\right)^{1/s}
 \qquad(n>n_0).
\]
This proves~\eqref{eq:horizontal-tail}. The same upper
$\ell_s$ estimate holds for the cosets $Q\phi_\gamma^0$,
because every finite initial sum belongs to $X_0$.
The uncountability and norm-span properties of the family
belong to the structural results whose proofs are omitted here;
see Remark~\ref{rem:scope-of-note}.
\end{remark}

\begin{proposition}[Bounded completeness of the basis of $E$]
\label{prop:combined-bounded-completeness}
Under the parameter growth condition~\eqref{eq:parameter-growth},
the transfinite basis $(e_\alpha)_{\alpha<\omega_1}$ of $E$
is boundedly complete.
\end{proposition}

\begin{proof}
Suppose, to the contrary, that there are strictly increasing
ordinals $(\alpha_i)_i$ and scalars $(a_i)_i$ for which
\[
 S_N=\sum_{i=1}^N a_i e_{\alpha_i},\qquad
 M:=\sup_N\|S_N\|_{2,s}<\infty,
\]
but $(S_N)_N$ does not converge in norm. Choose $\varepsilon>0$
and $0=q_0<q_1<q_2<\cdots$ such that
\[
 x_n=\sum_{i=q_{n-1}+1}^{q_n}a_i e_{\alpha_i},\qquad
 \|x_n\|_{2,s}\ge\varepsilon\quad(n\in\N).
\]
Put $P_m=\pi_{\{1,\ldots,m\}}$, $P_0=0$, and
\[
 D_n=[\alpha_{q_{n-1}+1},\alpha_{q_n}].
\]
The ordinal intervals $D_n$ are successive. We distinguish
two cases.

\smallskip
\noindent\emph{Case 1.} There are $m\in\N$, $\delta>0$, and
infinitely many $n$ for which $\|P_mx_n\|_{2,s}\ge\delta$.
Since
\[
 \|P_mx_n\|_{2,s}\le\sum_{k=1}^m\|\pi_kx_n\|_{2,s},
\]
there are a fixed $k_0\le m$ and a subsequence $(x_{\nu_i})_i$
such that, with $\eta=\delta/m$,
\[
 z_i=\pi_{k_0}x_{\nu_i},\qquad
 \|z_i\|_{E_{k_0}}=\|z_i\|_{2,s}\ge\eta.
\]
Here we used Proposition~\ref{prop:combined-norm}.
Choose $g_i\in G_{k_0}$ supported in $D_{\nu_i}\cap C_{k_0}$
with $g_i(z_i)>\eta/2$. Such a choice is possible by symmetry,
the norming property, and closure under ordinal interval restrictions.
The functionals $g_i$ are successive in the component order.

Put $p=p_{k_0}$ and $p'=p/(p-1)$. For every $j$, choose
a rational $t_j$ with
\[
 \tfrac12 n_j^{-1/p}\le t_j\le n_j^{-1/p}.
\]
The component $(p,j)$-operation gives
\[
 g=\frac{t_j}{m_j}\sum_{i=1}^{n_j}g_i\in G_{k_0}\subset G.
\]
Set $N_j=q_{\nu_{n_j}}$. The ordinal support restrictions give
$g_i(S_{N_j})=g_i(z_i)$ for $i\le n_j$. Therefore
\[
 M\ge g(S_{N_j})
 =\frac{t_j}{m_j}\sum_{i=1}^{n_j}g_i(z_i)
 >\frac{\eta}{4}\,\frac{n_j^{1/p'}}{m_j}.
\]
This contradicts~\eqref{eq:parameter-growth}.

\smallskip
\noindent\emph{Case 2.} For every $m\in\N$,
$\|P_mx_n\|_{2,s}\to0$ as $n\to\infty$.
Inductively choose a subsequence $(x_{\nu_i})_i$ and integers
$0=b_0<b_1<b_2<\cdots$ so that
\[
 \|P_{b_{i-1}}x_{\nu_i}\|_{2,s}
 <\varepsilon 2^{-i-2},\qquad
 \operatorname{ind}(x_{\nu_i})\subset\{1,\ldots,b_i\}.
\]
The second choice is possible because $x_{\nu_i}$ has finite support.
Set
\[
 J_i=\{b_{i-1}+1,\ldots,b_i\},\qquad
 y_i=\pi_{J_i}x_{\nu_i}.
\]
Then $(y_i)_i$ is an h-block sequence and
\[
 \|x_{\nu_i}-y_i\|_{2,s}<\varepsilon 2^{-i-2},\qquad
 \|y_i\|_{2,s}>\varepsilon/2.
\]
Choose $h_i\in G$ with $h_i(y_i)>\varepsilon/4$, and restrict
$h_i$ to $D_{\nu_i}\cap C_{J_i}$. By the two restriction rules,
the resulting $f_i$ belongs to $G$ and has the same value on $y_i$.
The functionals $(f_i)_i$ form an h-block sequence.

For every $j$, choose a rational $u_j$ with
\[
 \tfrac12 n_j^{-1/2}\le u_j\le n_j^{-1/2}.
\]
The horizontal $(2,j)$-operation gives
\[
 f=\frac{u_j}{m_j}\sum_{i=1}^{n_j}f_i\in G.
\]
Set $N_j=q_{\nu_{n_j}}$. The support restrictions imply
$f_i(S_{N_j})=f_i(y_i)$ for $i\le n_j$, so
\[
 M\ge f(S_{N_j})
 =\frac{u_j}{m_j}\sum_{i=1}^{n_j}f_i(y_i)
 >\frac{\varepsilon}{8}\,\frac{n_j^{1/2}}{m_j}.
\]
This contradicts~\eqref{eq:parameter-growth} with $q=2$.
Both cases are impossible, proving bounded completeness.
\end{proof}

Since $G\subset B_{X_0}$ and $G$ norms $E$, the space $X_0$
is $1$-norming for $E$. Proposition~\ref{prop:bounded-complete-predual}
now gives the canonical onto isometry
\[
 E\longrightarrow X_0^*,\qquad
 x\longmapsto\bigl(h\mapsto h(x)\bigr).
\]
Thus property~\ref{combined:basis} follows from the norming-set
construction and the preceding proposition.

\begin{theorem}[Summary of the properties of $E$]
\label{thm:combined-summary}
Let $E$ be the space defined by $G$, with the parameters and
partition fixed in Section~\ref{sec:norming-combined}.
Let $\Phi^0=\{\phi_\gamma^0:\gamma\in\Gamma_0\}$ be the
family of infinite horizontal special functionals outside $X_0$.
For $k\ge1$, let
$\Phi^k=\{\phi_\gamma^k:\gamma\in\Gamma_k\}$ be the
family of infinite special functionals of $E_k$, extended to
$E$ by $\iota_k$. We call the families $\Phi^k$, $k\ge1$,
the \emph{vertical special families}. Put $Q:E^*\to E^*/X_0$.
Then the following hold.
\begin{enumerate}[label=\textup{(\alph*)},leftmargin=*]
\item\label{summary-E:upper}
The basis $(e_\alpha)_{\alpha<\omega_1}$ admits an upper
$\ell_2$ estimate with constant one: for every finitely
supported scalar family $(a_\alpha)$,
\[
 \left\|\sum_\alpha a_\alpha e_\alpha\right\|_{2,s}
 \le\left(\sum_\alpha|a_\alpha|^2\right)^{1/2}.
\]
\item\label{summary-E:hblocks}
For every normalized infinite h-block sequence $(x_n)_n$,
the space $\clspan\{x_n:n\in\N\}$ has nonseparable dual.
\item\label{summary-E:span}
The biorthogonal functionals, together with the horizontal and
vertical special families, generate the dual in norm:
\[
 E^*
 =\clspan\left(\{e_\alpha^*:\alpha<\omega_1\}
                      \cup\bigcup_{k\ge0}\Phi^k\right)
 =\clspan\left(\bigcup_{k\ge1}E_k^*\cup\Phi^0\right).
\]
\item\label{summary-E:quotient}
The horizontal quotient family $(Q\phi_\gamma^0)_\gamma$
admits an upper $\ell_s$ estimate, and for every $k\ge1$
the vertical quotient family $(Q\phi_\gamma^k)_\gamma$
admits an upper $\ell_{r_k}$ estimate, all with constant one.
Equivalently, setting $q_0=s$ and $q_k=r_k$ for $k\ge1$,
\[
 \left\|\sum_{\gamma\in F}a_\gamma Q\phi_\gamma^k\right\|
 \le\left(\sum_{\gamma\in F}|a_\gamma|^{q_k}\right)^{1/q_k}
\]
for every $k\ge0$, every finite $F\subset\Gamma_k$, and
every scalar family $(a_\gamma)_{\gamma\in F}$. Moreover,
\[
 E^*/X_0
 =\clspan\{Q\phi_\gamma^k:k\ge0,\ \gamma\in\Gamma_k\}.
\]
\end{enumerate}
\end{theorem}

\begin{proof}[Proof of the estimates]
For~\ref{summary-E:upper}, we first check that
$\|f\|_{\ell_2(\omega_1)}\le1$ for every $f\in G$.
At the initial stage, $f\in G_k$ for some $k$, and
Lemma~\ref{lem:norming-pr-lp}, transported to $C_k$, gives
\[
 \|f\|_{\ell_2}\le\|f\|_{\ell_{p_k}}\le1,
\]
since $p_k<2$. For a horizontal $(2,j)$-operation, disjointness
of the children's supports gives
\[
 \left\|\frac1{m_j}\sum_t\mu_t f_t\right\|_2^2
 =\frac1{m_j^2}\sum_t|\mu_t|^2\|f_t\|_2^2
 \le\frac1{m_j^2},
\]
provided the bound holds for the children. Thus a horizontally
weighted functional has $\ell_2$ norm at most the reciprocal
of its chosen weight. The special-sum rules preserve the unit
bound by the triangle inequality and $\sum_j1/m_j<1$;
in rule~\textup{(5)}, the retained weight sets are disjoint
and each coefficient has absolute value at most one.
Symmetry, both types of restriction, and rational convex
combinations also preserve the bound. Induction as in
Remark~\ref{rem:combined-induction} proves the claim.
Applying the Cauchy--Schwarz inequality and taking the supremum
over $f\in G$ gives~\ref{summary-E:upper}.

The horizontal estimate in~\ref{summary-E:quotient} was proved in
Remark~\ref{rem:infinite-horizontal-functionals}. The vertical
estimates follow from Theorem~\ref{thm:pr-summary}\textup{(d)}
and the isometric embeddings $j_k:W_k\to W$ described in
Section~\ref{sec:operator-B}. The norm-span assertion in the
quotient follows from~\ref{summary-E:span}.

Assertions~\ref{summary-E:hblocks} and~\ref{summary-E:span}
require the corresponding HI and James-tree arguments for the
combined construction. The discussion and references following
Theorem~\ref{thm:pr-summary} apply here as well; for horizontal
constructions involving families of spaces, see
also~\cite{AF00,AR08}. These lengthy technical proofs are
omitted from the present paper.
\end{proof}

\begin{corollary}\label{cor:E-hereditary-nonseparable-dual}
Every infinite-dimensional closed subspace of $E$ has
nonseparable dual.
\end{corollary}

\begin{proof}
Let $Y\subset E$ be an infinite-dimensional closed subspace.
By~\cite[Proposition~1.3(3)]{ALT05}, $Y$ contains a closed
subspace isomorphic to a block subspace
$U=\clspan\{u_n:n\in\N\}$, where $(u_n)_n$ is normalized.
It suffices to show that $U^*$ is nonseparable: the same will
then hold for the dual of that subspace of $Y$, which is a
quotient of $Y^*$ by the Hahn--Banach theorem.

Suppose first that, for some $k$, the operator
$\pi_k|_U:U\to E_k$ is not strictly singular. There is then
an infinite-dimensional closed subspace $Z\subset U$ on which
$\pi_k$ is an isomorphism onto a closed subspace of $E_k$.
By Corollary~\ref{cor:pr-hereditary-nonseparable-dual}, applied
to the isometric copy $E_k=E_{p_k,r_k}(C_k)$, this image has
nonseparable dual. Hence $Z^*$ is nonseparable, and the
surjective restriction map $U^*\to Z^*$ shows that $U^*$ is
nonseparable.

It remains to consider the case where every $\pi_k|_U$ is
strictly singular. Put $P_m=\sum_{k=1}^m\pi_k$ and $P_0=0$.
Finite sums of strictly singular operators are strictly singular,
so each $P_m|_U$ is strictly singular. Choose positive numbers
$(\varepsilon_n)_n$ with
$\theta:=\sum_n\varepsilon_n<1/2$.
We inductively choose successive normalized finite blocks
$(v_n)_n$ of $(u_n)_n$ and integers
$0=b_0<b_1<b_2<\cdots$ such that
\[
 \|P_{b_{n-1}}v_n\|_{2,s}<\varepsilon_n,
 \qquad \operatorname{ind}(v_n)\subset\{1,\ldots,b_n\}.
\]
To choose $v_n$, restrict $P_{b_{n-1}}$ to the
infinite-dimensional tail of $U$ following the preceding block.
Strict singularity gives a unit vector whose image has arbitrarily
small norm. Approximation by a finite block, followed by
normalization, gives the required $v_n$. Then choose $b_n$
to contain its finite horizontal support.

Set
\[
 w_n=\pi_{\{b_{n-1}+1,\ldots,b_n\}}v_n.
\]
The vectors $w_n$ are successive in the ordinal order, their
horizontal supports are successive, and
\[
 \|v_n-w_n\|_{2,s}<\varepsilon_n,
 \qquad \|w_n\|_{2,s}>1-\varepsilon_n>1/2.
\]
Thus $(w_n)_n$ is a seminormalized h-block sequence.
To verify explicitly the equivalence of $(v_n)_n$ and $(w_n)_n$,
let $(a_n)_n$ be finitely supported. Bimonotonicity and
$\|v_n\|_{2,s}=1$ imply
\[
 \sup_n|a_n|\le\left\|\sum_n a_nv_n\right\|_{2,s}.
\]
Consequently,
\[
 \left\|\sum_n a_n(v_n-w_n)\right\|_{2,s}
 \le\sum_n|a_n|\varepsilon_n
 \le\theta\left\|\sum_n a_nv_n\right\|_{2,s}.
\]
Since $\theta<1$, the map $v_n\mapsto w_n$ extends to an
isomorphism between their closed linear spans. After normalizing
$(w_n)_n$, Theorem~\ref{thm:combined-summary}\textup{(b)} shows
that its closed span has nonseparable dual. The same is true
of $\clspan\{v_n:n\in\N\}\subset U$. Hahn--Banach again
implies that $U^*$ is nonseparable, completing the proof.
\end{proof}

\begin{remark}[Scope of the note]\label{rem:scope-of-note}
The norming sets $G_{p,r}$ and $G$ have been defined explicitly,
so the definitions of the spaces $E_{p,r}$ and $E$ are complete.
We have proved directly the basic properties of their bases,
including bounded completeness and the canonical predual
identifications, the isometric identification of the components
of $E$, and the estimates for special functionals and their tails.
Starting from the structural properties of $E$ stated in
Section~\ref{sec:construction}, we have also given the full
construction of $B$ under CH and the proof that
$X=\ker(Q-B)$ has no infinite-dimensional separable quotient.

The construction results whose proofs are omitted here are the
block-subspace and dual norm-span assertions in
Theorems~\ref{thm:pr-summary}\textup{(b),(c)}
and~\ref{thm:combined-summary}\textup{(b),(c)}, together with
the absence of isomorphic copies of $\ell_1$ in these spaces.
Corollaries~\ref{cor:pr-hereditary-nonseparable-dual}
and~\ref{cor:E-hereditary-nonseparable-dual} prove the passage
from the block-subspace assertions to nonseparable duals of all
infinite-dimensional closed subspaces. The dual norm-span
assertions yield the density of the linear span of the special
cosets in the corresponding quotients. For $E$, this means
\[
 E^*/X_0
 =\clspan\{Q\phi_\gamma^k:
                 k\ge0,\ \gamma\in\Gamma_k\},
\]
with the analogous statement for each component space.
These structural properties are used as inputs in the operator
argument; their derivation from the norming sets is not included
in this note.

The omitted proofs require technical methods from the theory of
hereditarily indecomposable Banach spaces. Their inclusion would
take us beyond the purpose of this paper, which is to make the
construction of the space $X$ and the mechanism behind its
failure to have a separable quotient more accessible.
If necessary, these proofs will be presented in a separate paper.
\end{remark}

\end{document}